\documentclass[a4paper]{article}
\usepackage{amssymb, amsmath, enumerate, geometry,epsf}
\usepackage[all]{xy}
\usepackage{amsfonts}
\usepackage{latexsym}
\usepackage{amssymb}
\usepackage{amscd}

\newtheorem{thm}{Theorem}[section]
\newtheorem{defn}[thm]{Definition}
\newtheorem{dfn}[thm]{Definition}

\newtheorem{lemma}[thm]{Lemma}
\newtheorem{prop}[thm]{Proposition}
\newtheorem{cor}[thm]{Corollary}

\newtheorem{example}[thm]{Example}
\newtheorem{remark}[thm]{Remark}

\newenvironment{pf}{\emph{Proof.}}{\hfill$\Box$\\[1mm]}

\def\a{\alpha} \def\b{\beta}  \def\d{\delta}
\def\e{\varepsilon} \def\f{\varphi}  \def\l{\lambda}
 \def\r{\rho} \def\s{\sigma}  
\def\xip{\xi'}
   \def\D{\Delta}

\def\ss{\mathbf{s}}
\def\f{\mathbf{f}} 
\def\q{\mathbf{q}} 
 \def\l{\lambda}

\def\r{\mathbf{r}}   

\def\bC{\mathbb{C}}
\def\bN{\mathbb{N}}

\def\bZ{\mathbb{Z}}

\def\one{\mathbf{1}}
\def\bigone{\mathbf{1}}

\def\A{{A}}

\def\B{B}
\def\C{\mathcal{C}}

 \def\H{H}

\def\M{M}

\def\P{\Psi}

\def\RM{R(M)}

\def\Ad{\mathrm{Ad}} \def\id{\mathrm{id}}
 \def\tr{\mathrm{tr}}

\def\im{\mathrm{im}}
\def\ker{\mathrm{ker}}
\def\coker{\mathrm{coker}}
\def\Ext{\mathrm{Ext}}
\def\Tor{\mathrm{Tor}}
\def\Hom{\mathrm{Hom}}

\def\otherwise{\mathrm{otherwise}}

\def\op{\mathrm{op}}
\def\cop{\mathrm{cop}}
\def\Ob{\mathrm{Ob}}

\def\a{\alpha}
\def\b{\beta}
\def\bp{b'}
\def\slq2{SL_q (2)}
\def\s{\sigma}

\def\to{\rightarrow}

\def\tr{\triangleright}
\def\tl{\triangleleft}
\def\btr{\blacktriangleright}
\def\btl{\blacktriangleleft}

\def\Mor{\mathrm{Mor}}

\def\kchoosel{\left(\begin{array}{cr} k\\ l \end{array} \right)}

\def\span{\mathrm{span}}
\def\Ext{\mathrm{Ext}}

\def\zero{{{(0)}}}
\def\one{{{(1)}}}
\def\two{{{(2)}}}
\def\three{{{(3)}}}
\def\four{{{(4)}}}

\def\tr{\triangleright}
\def\tl{\triangleleft}
\def\btr{\blacktriangleright}
\def\btl{\blacktriangleleft}

\def\atr{\triangleright}
\def\atl{\triangleleft}
\def\abtr{\blacktriangleright}
\def\abtl{\blacktriangleleft}

\def\abtr{\triangleright}
\def\abtl{\triangleleft}
\def\atr{\blacktriangleright}
\def\atl{\blacktriangleleft}

\def\field{\bC}
\def\kVec{\bC\!\!-\!\!{\bf\mathrm{Vec}}}

\begin{document}

\title{Braided  homology of quantum groups}
\author{Tom~Hadfield${}^1$, Ulrich~Kr\"{a}hmer${}^2$} 
\date{\today}
\maketitle

\centerline{${}^1$ Thomas.Daniel.Hadfield@gmail.com} 
\centerline{Supported  by EU Transfer of Knowledge contract MKTD-CT-2004-509794}
\centerline{and by an EPSRC Postdoctoral Fellowship}
\centerline{}
\centerline{ ${}^2$ University of Glasgow, Department
of Mathematics}
\centerline{University Gardens, G12 8QW Glasgow, UK, 
ukraehmer@maths.gla.ac.uk}
\centerline{Supported by an 
EU Marie Curie and an EPSRC Postdoctoral Fellowship.}
\centerline{}
\centerline{MSC (2000): 58B32, 16W30, 19D55} 
\centerline{}
\centerline{Keywords: Hochschild homology, cyclic homology, braided homology, quantum groups}

\begin{abstract}
We study braided Hochschild and cyclic homology of ribbon algebras in braided monoidal categories, 
 as introduced by Baez and by Akrami and Majid.
 We compute this invariant for several examples coming from quantum groups and braided groups. 
\end{abstract}

\section{Introduction}
Braided Hochschild homology $HH^\Psi_\ast (\A)$
 was first defined by Baez \cite{baez} for an algebra $\A$ in a braided monoidal category $\C$.  
 Baez defined $HH^\Psi_\ast (\A)$  by an explicit chain complex, and showed that $HH^\Psi_\ast (A) = \Tor^{\A \hat{\otimes} \A^\op}_\ast (\A,\A)$, where  $\A \hat{\otimes} \A^\op$ is the braided enveloping algebra of $\A$.
  However, his constructions relied on the assumption that 
 $\A$ is weakly $\Psi$-commutative, meaning that ${\mu} \circ \Psi^2 = {\mu}$, where ${\mu}$ is the multiplication morphism of $\A$.
  This was overcome  by Akrami and Majid \cite{am}, who replaced weakly $\Psi$-commutative algebras by ribbon algebras,  namely algebras $\A$ in  $\C$ possessing 
   an invertible morphism $\s \in \Mor_\C (\A,\A)$ that satisfies $\s \circ {\mu} = {\mu} \circ (\s \otimes \s) \circ \Psi^2$. 
   They showed that there is a corresponding cyclic theory (in the sense of Connes \cite{aCo83}), which they called braided cyclic homology $HC^{\Psi,\s}_\ast (\A)$.
   If $\C$ is the category of $\field$-vector spaces with braiding given by the flip  $v \otimes w \mapsto w \otimes v$, then an algebra  in $\C$ is simply an associative $\field$-algebra, ribbon automorphisms are precisely algebra automorphisms, and the Akrami-Majid cyclic object reduces to the 
   cyclic object defining twisted cyclic homology as studied by Kustermans, Murphy and Tuset \cite{kmt}.

In this paper we first show that the ribbon automorphism $\s$ is also the missing ingredient required to make Baez's realisation of braided Hochschild homology as a derived functor work in full generality. 
We prove a braided analogue of a result of Feng and Tsygan \cite{FT},  
 namely  for Hopf algebras, Hochschild homology can be realised as a derived functor in the category of modules over the algebra itself. 
 We apply  this machinery to concrete examples of algebras and Hopf algebras in braided monoidal categories: the braided line, braided plane  and braided quantum $SL(2)$.

 From the viewpoint of braided monoidal categories, braided Hochschild and cyclic homology are the natural abstractions of the original definitions of Hochschild and of Connes, and therefore natural objects to  study. 
 The computations carried out in this paper also show that braided Hochschild homology is (similarly to twisted Hochschild homology \cite{bz, tompodles, hk, hk2, kmt, sitarz}) often less degenerate than classical Hochschild homology, for example in the sense that it overcomes the so-called dimension drop observed for quantisations of Poisson varieties \cite{FT}. 
 To show whether the standard applications of cyclic homology in noncommutative geometry admit generalisations to the braided setting seems a promising direction for future research.
 Another interesting question is the relevance of these
 invariants for the study of conformal field theories that can be
 described in terms of ribbon algebras, see 
e.g.~\cite{juergen}.
 
  A summary of this paper is as follows. 
  Throughout   we work over $\bC$ as ground field.
  In Section \ref{section:hom_alg_braided_cat} we recall the definitions of a ribbon algebra and a Hopf algebra in a braided monoidal category $\C$, with particular reference to the motivating example $\C = \C(\H)$, the category of comodules of  a coquasitriangular Hopf algebra $\H$. 
  In Section \ref{section:braided_hoch_and_cyclic} we define braided Hochschild homology $H^{\P,\s}_\ast (\A, \M)$ for  a ribbon algebra $(\A,\s)$ and an $\A$-bimodule $\M$ in  $\C$ as the homology of a specific complex.
   We show that for $\C = \C(\H)$ this can be realised as a derived functor  over an appropriate  braided enveloping algebra $\A^e$ (Theorem \ref{braided_Hoch_is_Tor}), thus generalising Baez's result to the case $\s \neq \id$. 
  We discuss the precise relation to  the  Akrami-Majid cyclic object associated to $(\A,\s)$.
    
  In Section \ref{section:HH_Hopf} we prove a braided analogue of a result of Feng and Tsygan \cite{FT}, namely that for  a Hopf algebra $\A$ in $\C(\H)$ there is an  isomorphism of vector spaces
$H_n^{\P, \s} (\A,\M) \simeq \Tor_n^\A (\RM,\bC)$, for a suitable right $\A$-module $R(\M)$ associated to any  $\A$-bimodule $\M$ (Theorem \ref{braided_ft}).  In Section \ref{section:braided_line} we apply this machinery to the  braided line and  braided plane.

In the final Section we consider braided Hopf algebras associated to coquasitriangular Hopf algebras  via the process known as transmutation.
 For quantum $SL(2)$, we obtain a no dimension drop 
type result (along the lines of \cite{bz,hk,hk2}) for the associated braided Hopf algebra $B$, namely that $H_n^{\P,\s} (B,B) =0$ for $n > 3$, and $H_3^{\P,\s} (B,B) \cong \bC$.

\section{Braided monoidal categories}
\label{section:hom_alg_braided_cat}

\subsection{Ribbon algebras and Hopf algebras in braided monoidal categories}
Recall \cite{shahn_book} that a monoidal category is a category $\C$ together with a functor 
$\otimes : \C \times \C \to \C$, an object $\bigone \in \Ob(\C)$ and isomorphisms of functors 
$\Phi : ( \cdot  \otimes \cdot ) \otimes \cdot \cong \cdot \otimes  ( \cdot \otimes \cdot )$ (the associator), $\ell :  \cdot \otimes \bigone \to \id$ and $r : \bigone \otimes \cdot \to \id$, assumed to satisfy certain consistency relations, which abstract the properties of the tensor product, say of vector spaces, over a field.
A monoidal category is  Ab-monoidal if for any $V$, $W \in \Ob(\C)$, the set $\Hom(V,W)$ is an additive abelian group such that composition and tensor product of morphisms are bilinear. 

\begin{figure}
\label{fig:vwbraid}
 \[ \epsfbox{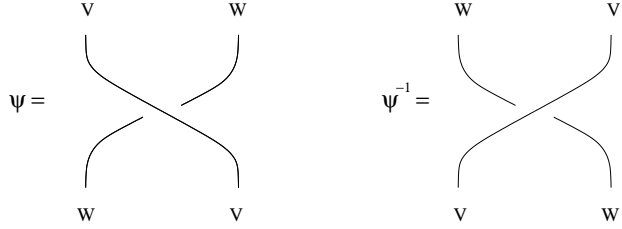} \]
\caption{The braiding $\P : V \otimes W \to W \otimes V$ and its inverse.}
\label{vwbraid}
\end{figure}

Given a monoidal category $\C$, we define a new functor $\otimes^\op : \C \times \C \to \C$ by $\otimes^\op (V,W) = W \otimes V$. 
 A braided monoidal category is a monoidal category $\C$ equipped with a braiding, an isomorphism of functors $\Psi : \otimes \to \otimes^\op$, 
  obeying the so-called hexagon relation. 
  As in \cite{am, baez}, we shall use standard graphical notation to depict morphisms in $\C$. 
   For example the braiding    $\P : V \otimes W \to W \otimes V$ and its inverse are shown in Figure  \ref{vwbraid}.
    Note that $\C$ equipped with $\Psi^{-1}$ is also a braided monoidal category. The choice of $\Psi$ rather than $\P^{-1}$ is simply a matter of convention. 

An algebra (or monoid) in $\C$  is an object $\A$ of $\C$ with  morphisms ${\mu} : \A \otimes \A \to \A$ (multiplication) and 
$\eta : \bigone \to \A$, such that  ${\mu} \circ ({\mu} \otimes \id) ={\mu} \circ (\id \otimes {\mu})$ 
and ${\mu} \circ (\eta \otimes \id) = \id = {\mu} \circ (\id \otimes \eta)$. 
Here and in the sequel we suppress $\Phi$, $\ell$, $r$. 
 Dually, a coalgebra (or comonoid) in $\C$ is an object $\A$ of $\C$ with morphisms $\D : \A \to \A \otimes \A$ (comultiplication) and $\e : \A \to \bigone$ (counit) satisfying $(\D \otimes \id) \circ \D = (\id \otimes \D) \circ \D$ and $(\e \otimes \id) \circ \D = \id = (\id \otimes \e) \circ \D$.
  Given an algebra $\A$ in $\C$, we define the opposite algebra $\A^\op$ to be the object $\A$ equipped with multiplication $\mu^\op = \mu \circ \P$. Similarly we define the coopposite coalgebra $\A^\cop$.

Given algebras $\A$ and $\B$ in $\C$, the braided tensor product algebra $\A \hat{\otimes} \B$ is defined to be the object $\A \otimes \B$  with multiplication 
$${\mu}_{\A \hat{\otimes} \B} = ({\mu}_\A \otimes {\mu}_\B) \circ ( \id \otimes \Psi  \otimes \id) : (\A \otimes \B) \otimes  (\A \otimes \B) \to \A \otimes \B$$

\begin{figure}
\label{fig:braided_coproduct}
\[ \epsfbox{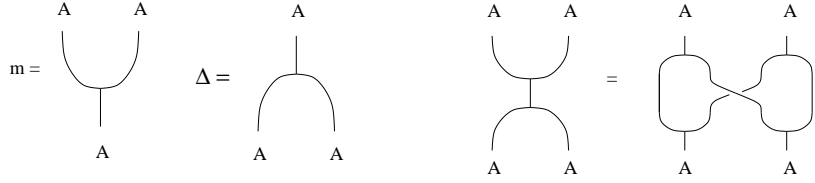} \]
\caption{Multiplication, comultiplication and the bialgebra condition}
\label{braided_coproduct}
\end{figure}

A bialgebra $\A$ in $\C$ is an algebra and coalgebra for which $\D : \A \to \A \hat{\otimes} \A$ and  $\e : \A \to \bigone$ are algebra morphisms (Figure \ref{braided_coproduct}), i.e. 
$\D \circ {\mu} = ({\mu} \otimes {\mu}) \circ (\id \otimes \Psi \otimes \id) \circ (\D \otimes \D)$, and $\e \circ {\mu} = \e \otimes \e$. 
Further, a   Hopf algebra (or braided group) in  $\C$ is a bialgebra $\A$ together with a 
 morphism $S : \A \to \A$ (antipode), satisfying ${\mu} \circ ( S \otimes \id) \circ \D = \eta \circ \e = {\mu} \circ (\id \otimes S) \circ \D$.
It follows that $S$ is a braided antihomomorphism, i.e. $S \circ {\mu} = {\mu} \circ \Psi \circ (S \otimes S)$.
 
Finally we need the notion of a ribbon algebra in a braided monoidal category, which is the basic data needed for Akrami and Majid's construction of braided cyclic (co)homology. 

\begin{figure}
\label{fig:ribbon}
 \[ \epsfbox{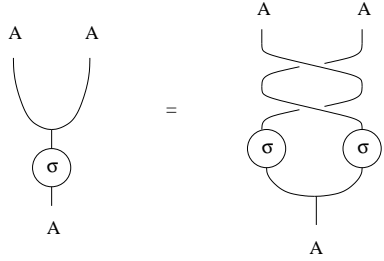} \]
\caption{The ribbon property $\s \circ {\mu} = {\mu} \circ (\s \otimes \s) \circ \P^2$.}
\label{ribbon}
\end{figure}

\begin{defn}\label{defn_ribbon_alg} (Figure \ref{ribbon}) A ribbon algebra in $\C$ is an algebra $(\A, {\mu},\eta)$ together with an invertible morphism $\s : \A \to \A$ (the ribbon automorphism) such that
\begin{equation}
\label{ribbon_aut}
{\mu} \circ (\s \otimes \s) \circ \P^2 = \s \circ {\mu} : \A \otimes \A \to \A, \quad \s\circ \eta = \eta
\end{equation}
and $(\s \otimes \id) \circ \P = \P \circ ( \id \otimes \s)$, $(\id \otimes \s) \circ \P = \P \circ (\s \otimes \id)$ on $\M \otimes \A$, $\A \otimes \M$ respectively, for all $\M \in \Ob (\C)$.
\end{defn}

Note that if $(\A,\s)$ is a ribbon algebra in $( \C, \P)$, then $(\A, \s^{-1})$ is a ribbon algebra in $(\C, \P^{-1})$.

\begin{example}\label{example_vs} If we take $\C = \kVec$ to be the category of complex vector spaces, then an algebra $\A$ in $\C$ is simply an ordinary algebra over $\field$. If  we take the braiding $\Psi$ to be the flip $v \otimes w \mapsto w \otimes v$,  then (\ref{ribbon_aut}) becomes $\s(ab) = \s(a) \s(b)$, i.e. $\s$ is an ordinary algebra automorphism of $\A$.
\end{example}

\subsection{Coquasitriangular Hopf algebras and braided monoidal categories}
\label{section_cqt}

We now recall the concepts of coquasitriangular and of
coribbon Hopf algebras. We refer to \cite{kassel,KS,shahn_book}
for more definitions and proofs. 
The most important  example for us is $H=\mathbb{C}_q[G]$ (defined in Example \ref{cqg}), the standard quantised coordinate ring of a complex simple Lie group $G$, for 
$q \in \mathbb{C} \setminus \{0\}$ not a root of unity.

\begin{dfn}
Let $H$ be a bialgebra in $\C = \kVec$. 
 Then $H$ is called coquasitriangular (cobraided) if there exist bilinear forms $\r$, $\bar \r$ on $\H$  such that for all $f,g,h \in H$ 
\begin{eqnarray}
&& \r(f_{(1)},g_{(1)}) \bar \r(f_{(2)},g_{(2)})=
	\bar\r(f_{(1)},g_{(1)}) \r(f_{(2)},g_{(2)})=
	\varepsilon(f)\varepsilon(g),\nonumber\\
&& gf=\r(f_{(1)},g_{(1)}) 
	f_{(2)}g_{(2)} \bar \r(f_{(3)},g_{(3)}),\nonumber\\
&& \r(fg,h)=\r(f,h_{(1)})\r(g,h_{(2)}),\quad
	\r(f,gh)=\r(f_{(1)},h)\r(f_{(2)},g).\nonumber
\end{eqnarray}
We call $\r$ a universal r-form on $\H$.
\end{dfn} 
 Here we use Sweedler's notation $\D(f) =  f_{(1)} \otimes f_{(2)}$ (summation suppressed) for the coproduct.
As a consequence $\r$ satisfies
$\r(1,f)=\r(f,1)=\varepsilon(f)$,  
and the quantum Yang-Baxter equation
 $\r_{12}\r_{13}\r_{23}=\r_{23}\r_{13}\r_{12}$. 
Here, and in the sequel, we use the convolution product of multilinear maps from a coalgebra to an algebra, and lower indices refer to the components in tensor products where these are applied. 
Thus explicitly the quantum Yang-Baxter relation reads
$$\r ( f_\one , g_\one ) \r( f_\two , h_\one) \r ( g_\two , h_\two ) = \r ( g_\one , h_\one ) \r( f_\one , h_\two ) \r ( f_\two , g_\two ) \quad \forall \; f,g,h \in H$$
The bilinear form
$\bar \r_{21}(f,g)=\bar \r(g,f)$ is also a 
universal r-form. If $H$ is  a Hopf algebra with antipode
$S$, then
$\bar \r(f,g)=\r(S(f),g)$, $\r(f,g)=\r(S(f),S(g))$.

Note that coquasitriangularity 
is the notion dual to the more commonly
used concept of quasitriangularity, which refers to the
existence of an element $R \in H \otimes H$ (the
``universal R-matrix'') with certain properties. While
in our case the multiplication and opposite
multiplication are linked via $\r$, a universal
R-matrix links the coproduct and the
coopposite one. The disadvantage of quasitriangularity
is that it often can only be applied rigorously when working
either with finite-dimensional Hopf algebras or with
topological ones and suitable completions of tensor
products. Working with the dual notion of
coquasitriangularity offers a   way out in many
examples. An alternative elegant and purely algebraic 
option is to pass to multiplier Hopf algebras, see for
example the recent work \cite{nt} by Neshveyev and
Tuset and the references therein.

For a bialgebra $\H$, let $\C(\H)$ denote the category of right $\H$-comodules, and $\C_f (\H)$ the category of finite-dimensional right $\H$-comodules. 
The algebra structure of $H$ corresponds to a monoidal
 structure on $\C(\H)$, $\C_f(H)$ with tensor product that  of underlying $\field$-vector spaces, and  coaction
$$V \otimes W \to V \otimes W \otimes \H, \quad v \otimes w \mapsto v_{\zero} \otimes w_{\zero} \otimes v_\one w_\one$$
The unit object is $\field$ 
 with trivial coaction
$\field \ni 1 \mapsto 1 \otimes 1 \in \field \otimes H$. 
A coquasitriangular 
structure on $H$ turns 
$\C(\H)$, $\C_f(H)$ into braided monoidal categories with braiding
\begin{equation}
\label{rightcomodbraiding}
	\Psi : V \otimes W \rightarrow W \otimes V,\quad
	v \otimes w \mapsto w_\zero \otimes v_\zero \; \r(v_\one,w_\one)
\end{equation} 
Then $\P^2 (v \otimes w) =  v_\zero \otimes w_\zero \q (v_\one,w_\one)$, where $\q = \r_{21} \r$.
 Replacing $\r$ by ${\bar{\r}}_{21}$ corresponds to replacing $\P$ by $\P^{-1}$. \

\begin{defn}
A coquasitriangular Hopf algebra $(H,\r)$ is called coribbon
if there exists an invertible central element 
$\ss \in H^\circ$ (the dual Hopf algebra) satisfying
$$
	\ss(fg)= \ss(f_{(1)}) \ss(g_{(1)}) \q(f_{(2)},g_{(2)}) ,\quad
	\ss(1)=1,\quad \ss(S(f))=\ss(f) \quad \forall \; f,\; g \in \H
$$
\end{defn}

Note that $\r(a, \cdot)$, $\r( \cdot, a) \in H^\circ$. 
Thus $\ss \q ( a, \cdot) = \q( a, \cdot) \ss$, and 
$\ss \q ( \cdot ,a) = \q( \cdot,a) \ss$, for all $a \in \H$.
For coribbon Hopf algebras, $\C_f(\H)$ becomes a ribbon category \cite{kassel} with ribbon structure 
$$\sigma : V \to V, \quad v \mapsto v_\zero \; \ss( v_\one ), \quad v \in \C_f(\H)$$
Since $\ss$ is cocentral, this is a morphism in $\C_f (\H)$. 
 It is also well-defined for objects of $\C(\H)$, and turns algebras in $\C(\H)$ into ribbon algebras in the sense of Definition \ref{defn_ribbon_alg}. 
 However  $\C(\H)$ is in general not a ribbon category, due to the lack of duals.
For a full discussion  see \cite{am}.

\begin{example} The trivial Hopf algebra $\H = \field$ is coquasitriangular, via $\r(1,1) =1$, and coribbon via $\ss(1) =1$. 
 Hence $\C(\H)$ is the braided monoidal category of Example \ref{example_vs}, however  the only ribbon automorphism arising from $\ss$ is $\s = \id$. 
\end{example}

\begin{example}\label{ex_canonical_braiding}
Any coribbon Hopf algebra $(\H,\r)$ becomes a  ribbon algebra in $\C(\H)$ with coaction $\Delta$. 
Then (\ref{ribbon_aut}) is a direct translation of the defining property of $\s$.
We call the resulting braiding 
\begin{equation}\label{canonical_braiding}
\P ( f \otimes g) = g_\one \otimes f_\one \r( f_\two , g_\two )
\end{equation}
 the canonical braiding on $(\H,\r)$.
\end{example}

\begin{example}\label{baez_example}
Any Hopf algebra $\H$ is a right $\H^\cop \otimes \H$-comodule algebra, via the coaction
$f \mapsto f_\two \otimes f_\one \otimes f_\three$.
 If $H$ is coquasitriangular, then $H^\cop \otimes H$ is also coquasitriangular, with universal r-form 
$( f \otimes g, h \otimes k ) \mapsto \bar{\r} (f,h) \r ( g , k)$. 
The resulting braiding on $H$ is 
\begin{equation}
\label{baez_braiding}
\P (f \otimes g) = \bar \r (f_\one , g_\one ) \; g_\two \otimes f_\two \; \r (f_\three , g_\three)
\end{equation}
By the definition of coquasitriangularity, we have
  ${\mu} \circ \P = {\mu}$. 
  So $\H$ is $\P$-commutative, in the terminology of \cite{baez}.
  In particular $\s = \id$ is a ribbon automorphism.
\end{example}

\begin{example}\label{cqg} 
Let $G$ be a complex simple Lie group, 
with Lie algebra $\mathfrak{g}$.
Let   $q \in \mathbb{C} \setminus \{0\}$ be not a 
root of unity, and 
$\H = \bC_q [G]$ and $U_q ( \mathfrak{g} )$ 
be the standard quantised coordinate ring 
and enveloping algebra of $G$ 
and $\mathfrak{g}$ respectively. 
Then $\C_{(f)} (\H)$ consists of (finite) direct 
sums of the so-called 
type I representations of $U_q ( \mathfrak{g})$ 
that are obtained by 
deformation of 
finite-dimensional representations of $G$. 
Both are braided monoidal categories. 
The corresponding universal r-form on $\bC_q [G]$ can
be given explicitly in terms of the so-called Rosso
form, see e.g. \cite{hodges}. 
For the classical matrix Lie groups $G$ the ribbon functionals on $\bC_q [G]$ were classified by Hayashi \cite{hayashi}.
  \end{example}
  
  \begin{example}\label{SLq2} Specialising Example \ref{cqg}, 
  the coordinate algebra $\H = \bC_q [ SL(2)] $ of quantum $SL(2)$  
  can be presented  in terms of generators $a$, $b$, $c$, $d$ with relations
\begin{equation}\label{slq2_one}
ab=qba, \quad
ac=qca, \quad
bd=qdb, \quad
cd=qdc, \quad
bc=cb, \quad ad-qbc=1, \quad
da - q^{-1} bc=1
\end{equation}
The Hopf algebra structure on $\H$ is given by
\begin{eqnarray}
&&
\Delta
\left[ \begin{array}{cc} a & b \cr c & d \end{array} \right]
= 
\left[ \begin{array}{cc} a & b \cr c & d \end{array} \right] 
\otimes
\left[ \begin{array}{cc} a & b \cr c & d \end{array} \right]
\nonumber\\
\label{slq2_two}
&&
S 
\left[ \begin{array}{cc} a & b \cr c & d \end{array} \right] = 
\left[ \begin{array}{cc} d & - q^{-1} b \cr -qc & a \end{array} \right],
\quad
\varepsilon 
\left[ \begin{array}{cc} a & b \cr c & d \end{array} \right] =
\left[ \begin{array}{cc} 1 & 0 \cr 0 & 1 \end{array} \right]
\end{eqnarray}

There are two useful 
vector space bases of $\H$: One is of
Poincar\'e-Birkhoff-Witt type and is given by the monomials 
$\{  e_{i,j,k}:=a^i b^j c^k,\quad i \in \bZ, \; j,k \in \bN  \}$,
where $a^i := d^{-i}$ for $i < 0$, and by convention
$x^0 =1$, for $x \in \H$, $x \neq 0$. The second relies
on the fact that $ \H$ is 
isomorphic as coalgebra to $ \mathbb{C} [SL(2)]$ and hence
cosemisimple (any comodule is the direct sum
of its irreducible subcomodules). Therefore, it
admits a Peter-Weyl type basis consisting of the matrix
coefficients $C^{(m)}_{rs}$, $m \in \mathbb{N}$,
$r,s=0,\ldots,m$, of the irreducible
corepresentations of $\H$. 
In particular,  
$\left[ \begin{array}{cc} C^{(1)}_{00} & C^{(1)}_{01} 
\cr C^{(1)}_{10} & C^{(1)}_{11} \end{array} \right] 
= \left[ \begin{array}{cc} a & b \cr c & d \end{array} \right]$.
The universal $r$-form $\r$ of Example \ref{cqg} is 
given on generators by 
\begin{equation}
\label{slq2_rform}
\left[\begin{array}{cccc}
\r(a , a) & \r( a,b) & \r(a ,c) & \r( a,d)\cr
\r(b , a) & \r(b ,b) & \r(b ,c) & \r( b,d)\cr
\r( c, a) & \r(c ,b) & \r(c ,c) & \r(c ,d)\cr
\r(d , a) & \r(d ,b) & \r(d ,c) & \r(d ,d)\cr
 \end{array}\right]
=
q^{-1/2} \left[\begin{array}{cccc}
 q & 0 & 0 & 1\cr
  0 & 0 & 0 & 0\cr
   0 & q - q^{-1}  & 0 & 0\cr
  1 & 0 & 0 & q \cr
 \end{array}\right]
\end{equation}
We refer to \cite{KS} for more details and proofs.
\end{example}

\subsection{Gauge transformations and cochain twists}\label{cochaintwist}

There is a natural notion of equivalence of braided monoidal categories, provided by  functors 
 implementing an equivalence of categories that transforms the braided monoidal structures into each other. 
We recall this explicitly for the case $\C = \C(\H)$, for $\H$ a coquasitriangular Hopf algebra. 
We refer to \cite{am,bm, kassel} for the precise definition in abstract braided monoidal categories.

 Let $H$ be a Hopf algebra, and $\f$ a convolution invertible bilinear form on $H$, with inverse $\bar \f$.
The multiplication of $H$ can be twisted by $\f$ as follows:
$$a \cdot_\f b = \f( a_\one , b_\one )\; a_\two  b_\two \; \bar \f ( a_\three, b_\three)$$
In  general the result is a nonassociative algebra $H^\f$ (a so-called quasi-Hopf algebra) since
$$( a \cdot_\f b ) \cdot_\f c = \partial \f ( a_\one, b_\one, c_\one ) \; a_\two \cdot_\f ( b_\two \cdot_\f c_\two ) \; \overline{\partial  \f} ( a_\three , b_\three, c_\three ),$$
where $\partial \f (a,b,c) = \f ( a_\one , b_\one) \f( a_\two b_\two , c_\one) \bar \f ( a_\three , b_\three c_\two ) \bar \f ( b_\four , c_\three)$. If $\partial \f = \e \otimes \e \otimes \e$, then $\f$ is said to be a 2-cocycle \cite{dt,drinfeld,shahn_book} and $H^\f$ is associative.
 We say that $\f$ is a 2-coboundary if $\f (x,y) = \varphi ( x_\one) \varphi (y_\one) \bar{\varphi} ( x_\two y_\two)$ for some linear map $\varphi : \H \to k$. 
  A 2-coboundary is automatically a 2-cocycle, but not conversely.
   For $\H^\f$ to be associative it is sufficient that $\partial \f$ is cocentral, meaning
   $$\partial \f ( a_\one , b_\one , c_\one ) \; ( a_\two \otimes b_\two \otimes c_\two ) = ( a_\one \otimes b_\one \otimes c_\one ) \; \partial \f ( a_\two , b_\two , c_\two )$$
   and not necessarily a 2-cocycle. 
 Since $H = H^\f$ as coalgebras,  $\C(H)$ and $\C( H^\f)$ coincide. 
 The two  algebra products correspond to different monoidal structures on this category.
 Both are given by  tensor product of vector spaces, but the twist is encoded in the isomorphism of vector spaces 
 $$F : V \otimes W \to V \otimes W, \quad  v \otimes w \mapsto v_\zero \otimes w_\zero \; \bar \f  ( v_\one , w_\one)$$
which transforms $V \otimes W$ as an object of $\C(H)$ into an object of $\C(H^\f)$. 
 Following \cite{drinfeld} we refer to $F$ as a  gauge transformation.
  This transforms an $H$-comodule algebra $A \in \C(H)$ to an $H^\f$-comodule algebra $A^\f \in \C( H^\f)$. $A$ and $A^\f$ coincide as vector spaces, but the multiplication in $A^\f$ is given in terms of the multiplication in $A$ by 
 ${\mu}^\f = {\mu} \circ F$. 
 For this  to be associative  $\f$ must be a 2-cocycle.
 Further, the r-form $\r$ on $\H$ transforms to $\f_{21} \r \bar{\f}$, which is an r-form on $\H^\f$ provided $\f$ is a 2-cocycle.
  For general $\f$ the formulae have to be modified by incorporating $\partial \f$, the Drinfeld associator. 
Twisting Hopf algebras and comodule algebras by coboundaries yields isomorphic structures:

\begin{lemma}\label{banal}
If $\f = \partial \varphi$, then $\H \cong \H^\f$ as Hopf algebras, and $\A \cong \A^\f$ as $\H$-comodule algebras.
\end{lemma}
\begin{pf} The isomorphisms are given by $h \mapsto \bar{\varphi} ( h_\one ) h_\two \varphi( h_\three)$, 
$a \mapsto a_\zero \varphi( a_\one )$
\end{pf}
Finally we note that ribbon functionals and automorphisms are preserved under gauge transformations.

\begin{example} Let $(\H,\r)$ be a coquasitriangular Hopf algebra, and $\q = \r_{21} \r$ as in Section \ref{section_cqt}. 
 Then an invertible central element $\ss \in \H^\circ$ is a ribbon functional if and only if $\q = \partial \bar{\ss}$. 
 In particular $\q$ is then a 2-cocycle. 
\end{example}

 Taking this into account, we get:

\begin{prop} Let $(\H,\r)$ be a coquasitriangular Hopf algebra such that $\r = \partial \varphi$. Suppose further that $\varphi^2$ is cocentral (i.e. in the centre of the dual Hopf algebra). 
Then ${\bar{\varphi}}^2$ is a ribbon functional. 
\end{prop}
\begin{pf}
Indeed, if $\r=\partial \varphi$,
then $\q=\partial \varphi^2$:
\begin{eqnarray}
	\q(a,b)
&=& \r(b_{(1)},a_{(1)})\r(a_{(2)},b_{(2)}) 
= \varphi (b_{(1)})\varphi(a_{(1)}) \bar\varphi(b_{(2)}a_{(2)})\r(a_{(3)},b_{(3)})\nonumber\\ 
&=& \varphi(b_{(1)})\varphi(a_{(1)}) \r(a_{(2)},b_{(2)}) \bar\varphi(a_{(3)}b_{(3)})
	\bar\r (a_{(4)},b_{(4)}) \r(a_{(5)},b_{(5)})\nonumber\\ 
&=&  \varphi (a_{(1)})\varphi(a_{(2)})\varphi(b_{(1)})\varphi(b_{(2)})\bar\varphi(a_{(3)}b_{(3)})
	\bar\varphi(a_{(4)}b_{(4)}) = \partial \varphi^2(a,b).\nonumber
\end{eqnarray}
\end{pf}

\begin{example}\label{rform_cobdy} 
The universal r-form of  $\bC_q [G]$ is a 2-coboundary. This is essentially shown in \cite{ksoib} Corollary 4.1.7, where the authors work with formal deformation quantisations. 

Concretely, one can prove directly that for 
$\bC_q [SL(2)]$ 
 the (convolution invertible) linear functional $\varphi$ defined by
\begin{equation}\label{varphi_SLq2}
 \varphi ( a^{i+1} b^j c^k ) = 0 = \varphi( d^{i+1} b^j c^k ), \; \varphi( b^j c^k ) = q^{(j+k)(1-j-k)/4} \b^j \gamma^k , \; i,j,k \geq 0, \; \beta \gamma = - q^{-3/2}
 \end{equation}
 satisfies $\r = \partial \varphi$, for the r-form
 $\r$ defined  by (\ref{slq2_rform}). 
For $ \beta = q^{-1/4},\gamma = -q^{-5/4}$ one obtains 
$ \varphi $ for which $ \varphi^2$ is cocentral, i.e., is the
 inverse of the ribbon functional. In terms of the
 Peter-Weyl basis, this $ \varphi $ is given by
$$
	\varphi (C^{(m)}_{rs})=\left\{
\begin{array}{ll}
(-1)^r q^{-r -{\tiny{\frac{1}{4}}} m^2 } \quad & r+s=m \\
0 \quad & \otherwise
\end{array}\right., 
$$ 
see \cite{chari} p.262. The inverse of the ribbon
 functional is thus given by
\begin{equation}\label{esistspaet}
	\varphi ^2 (C^{(m)}_{rs})=	(-1)^m q^{- {\tiny{\frac{1}{2}}} m^2-m}  
		  \delta_{rs}
\end{equation} 
Note that this is cocentral. The functional $ \varphi $
gives rise to the extension of the quantised universal
enveloping algebra dual to $ \mathbb{C}_q[SL(2)]$
known as the quantum Weyl group. Usually, the fact that
the coboundary of $ \varphi $ is $\r$ is
considered the other way round, the aim being to 
use the r-form to 
compute the coproduct of $ \varphi $ as in 
Proposition~8.2.3 in \cite{chari}. But as we see here,
the inverse point of view is relevant as well.  
\end{example}

\begin{example}\cite{diparo}
Let $G$ be a finite group and 
$\omega : G \times G \times G \rightarrow U(1)$
be a 3-cocycle (where $G$ acts trivially on $U(1)$). Define 
$$
		  \theta_g(x,y):=
		  \frac{\omega (g,x,y) \omega
		  (x,y,y^{-1}x^{-1}gxy)}
		  {\omega (x,x^{-1}gx,y)},\>
		  \gamma_x(g,h):=
		  \frac{\omega (g,h,x) \omega
		  (x,x^{-1}gx,x^{-1}hx)}
		  {\omega (g,x,x^{-1}hx)},\>
		  g,h,x,y \in G.
$$
Then the tensor product 
$A \otimes A^\circ$ of the
commutative algebra $A=\mathbb{C} [G]$ of
functions on $G$ with the group algebra
$A^\circ=\mathbb{C} G$ can be turned into a generalised
Drinfeld double with product 
$$
		  (\delta_g \otimes x)(\delta_h \otimes y):=
		  \delta_{g}(xhx^{-1}) \delta_g \otimes 
		  xy \theta_g(x,y)
$$ 
and coproduct
$$
		  \Delta (\delta_g \otimes x):=
		  \sum_{h,k \in G,hk=g} \delta_h \otimes x
		  \otimes \delta_k \otimes x \gamma_x(h,k), 
$$
where $x,y$ are identified with the corresponding basis
elements of $ \mathbb{C} G$, and $ \delta_g : G
 \rightarrow \mathbb{C}$ is the Dirac function.
This gives a quasi Hopf algebra $D^\omega (G)$
(the product is associative, but $ \Delta $ is only
quasi coassociative in general). Let 
$H:=(D^\omega (G))^\circ$ be the dual. Then, 
$H$ depends up to a 2-cochain twist only on the
cohomology class $ [\omega] \in H^3(G,U(1))$, 
see \cite{diparo}, and it is always coquasitriangular.
The comodule category describes
certain rational conformal field theories that arose in
the work of Dijkgraaf and Witten, and
this link to mathematical physics 
is another main motivation for the current
interest in ribbon categories, see e.g.~\cite{juergen}.    
\end{example}

\section{Braided Hochschild and cyclic homology}
\label{section:braided_hoch_and_cyclic}

\subsection{Braided Hochschild homology} 
\begin{figure}
\label{fig:dj}
 \[ \epsfbox{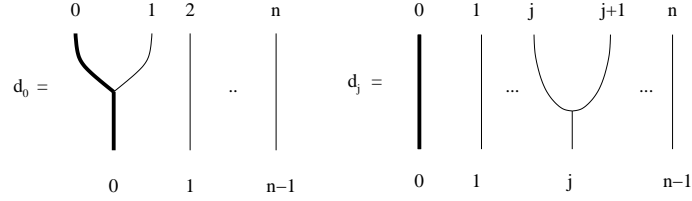} \]
\caption{The maps $d_0$ and $d_j$, $1 \leq j \leq n-1$}
\label{dj}
\end{figure}

\begin{figure}
\label{fig:dn}
 \[ \epsfbox{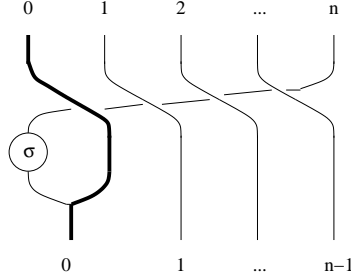} \]
\caption{The map $d_{n}$}
\label{dn}
\end{figure}

Let $(\A,\s)$ be a ribbon algebra in a braided monoidal Ab-category $\C$, which we assume possesses kernels and cokernels. 
This is in particular true for $\C = \C(\H)$, with $\H$ coquasitriangular.
There are obvious notions of left and right $\A$-module and $\A$-bimodule in $\C$. 
For example, by
  an $\A$-bimodule, we mean an object $\M$ of $\C$, together with morphisms (left and right actions)
  $\abtr \;:\; \A \otimes \M \to \M$, $\abtl \;:\; \M \otimes \A \to \M$
such that $\abtr  ( {\mu} \otimes \id ) = \abtr ( \id \otimes \abtr)$, $\abtl  ( \id \otimes {\mu}) = \abtl  ( \abtl \otimes \id)$,
$\abtr  (\id \otimes \abtl) = \abtl  (\abtr \otimes \id)$, 
and 
\begin{eqnarray}
&&\P  (\id \otimes \abtr) = (\abtr \otimes \id)  (\id \otimes \P )   ( \P \otimes \id) \quad \mbox{ on } \A \otimes \A \otimes \M
\nonumber\\
&&\P  (\id \otimes \abtl) = (\abtl \otimes \id)  (\id \otimes \P)  (\P \otimes \id)
 \quad \mbox{ on } \A \otimes \M \otimes \A
\nonumber
\end{eqnarray}
and similarly replacing $\P$ by $\P^{-1}$. 
From this point for convenience we drop the notation ``$\circ$" for composition of morphisms. 
Given such $\M$, define 
 $C_0 = C_0 (\A,\M)  = \M$,  $C_n = C_n (\A,\M) = \M \otimes \A^{\otimes n}$ for $n \geq 1$.
 In a diagram representing a morphism with source $\A^{\otimes m} \otimes \M \otimes \A^{\otimes n}$, we number the strands 0, 1,$\ldots$, $m+n$ (with $\M$ appearing as the strand labelled $m$, and represented graphically in bold).
We write
\begin{eqnarray}
&&\abtr_{m-1,m} = \id^{\otimes (m-1)} \otimes \abtr \otimes \id^{\otimes n} : \A^{\otimes m} \otimes \M \otimes \A^{\otimes n} \to \A^{\otimes (m-1)} \otimes \M \otimes \A^{\otimes n},\nonumber\\
&&\abtl_{m,m+1} = \id^{\otimes m} \otimes \abtl \otimes \id^{\otimes n} : \A^{\otimes m} \otimes \M \otimes \A^{\otimes n} \to \A^{\otimes m} \otimes \M \otimes \A^{\otimes (n-1)}\nonumber\\
&&{\mu}_{j,j+1} = \id^{\otimes j} \otimes {\mu} \otimes \id^{\otimes (m+n-j-1)}, \quad j \neq m-1,m\nonumber
\end{eqnarray}
For each $n \geq 1$, define maps $d_j : C_n \to C_{n-1}$, $0 \leq j \leq n$ by 
\begin{equation}
\label{d_j}
d_0 = \abtl_{\;0,1}, \quad d_j = {\mu}_{j,j+1} \quad 1 \leq j \leq n-1, \quad d_n = \abtr_{\;0,1} (\s \otimes \id^{\otimes n} ) \P_{[0,n-1],n}
\end{equation}
 where we  define, for $m \leq n < p$,  
\begin{eqnarray}
&&\P_{[m,n], n+1} = \P_{m,m+1} \P_{m+1,m+2} \ldots \P_{n,n+1},\nonumber\\
&&\P_{m,[m+1,n+1]} = \P_{n,n+1} \P_{n+1, n+2} \ldots \P_{m,m+1}\nonumber\\
&&\P_{[m,n],[n+1,p]} = \P_{[m+p-n-1,p-1],p} \ldots \P_{[m+1,n+1],n+2} \P_{[m,n], n+1}\nonumber\\
&&\P_{[m,n],[n+1,p]}^{-1} = (\P_{[m,n],[n+1,p]})^{-1}\nonumber
\end{eqnarray}
The maps $d_j$ are shown in Figures  \ref{dj} and  \ref{dn}.
 Together with maps $s_i : C_n \to C_{n+1}$, $0 \leq i \leq n$ defined by 
 $s_i = \id^{\otimes (i+1)} \otimes \eta \otimes \id^{\otimes (n-i)}$, this gives a simplicial object $\{ C_\ast (\A, \M) \}$. 
 Thus  $b_n : C_n \to C_{n-1}$ defined by  $b_n = \sum_{j=0}^{n} (-1)^j d_j$, gives a chain complex $\{ C_n , b_n \}_{n \geq 0}$.

\begin{defn} Braided Hochschild homology $H^{\P,\s}_\ast (\A,\M)$ of $\A$ with coefficients in $\M$  is  the homology of the complex $\{ C_n = \M \otimes \A^{\otimes n}, b \}_{n \geq 0}$.
\end{defn}
  
 \begin{example}\label{reduces_to_kmt} In the situation of Example \ref{example_vs}, a bimodule $\M$ over $\A$ is a bimodule in the usual sense, and  $H_n^{\P,\s} (\A, \M)$ reduces to $H_n ( \A, {}_\s \M)$, Hochschild homology of $A$ with coefficients in the $\A$-bimodule ${}_\s \M$, which is $\M$ as $\field$-module with bimodule structure $x \tr a \tl y = \s(x) \cdot a \cdot y$, where $\cdot$ is the original bimodule structure on $\M$.\\
 \end{example}
  
We now extend the derived functor interpretation of $H^{\P,\s}_\ast (\A,\M)$ from \cite{baez} to this  setting.

\begin{figure}
\label{fig:braidedenv}
 \[ \epsfbox{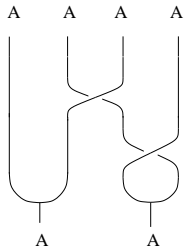} \]
\caption{The multiplication of the braided enveloping algebra $\A^e$.}
\label{braidedenv}
\end{figure}

\begin{defn}\label{defn_A^e}
The braided enveloping algebra $\A^e$  is the object $\A \otimes \A$ equipped with the multiplication morphism (Figure \ref{braidedenv})
\begin{equation}
\label{defn_btl_btr}
{\mu}_{\A^e} = ({\mu} \otimes {\mu}) \P^{-1}_{2,3} \P^{-1}_{1,2} : \A^{\otimes 4} \to \A^{\otimes 2}
\end{equation}
\end{defn}

That is, $\A^e$ is $\A \hat{\otimes} \A^\op$ taken in $( \C, \P^{-1})$ rather than $(\C, \P)$. 
We use this convention in order to be compatible with \cite{am}. 
Baez works with the opposite convention, taking $\A \hat{\otimes} \A^\op$ in $(\C, \P)$, resulting in a graphical calculus where braidings are replaced by inverse braidings.

\begin{figure}
\label{fig:leftrightaction}
 \[ \epsfbox{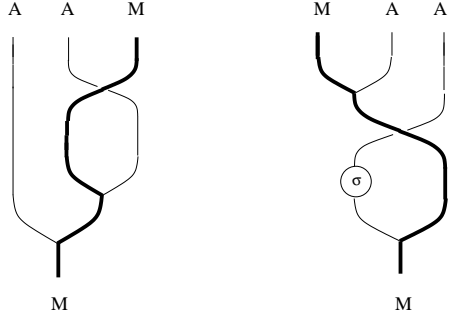} \]
\caption{The left and right actions of $\A^e$ on $\M$}
\label{leftrightaction}
\end{figure}

\begin{defn}\label{A_bimod_is_left_right_A^e_mod}
For any $\A$-bimodule $\M$, we define morphisms (Figure \ref{leftrightaction})
$$\atr = \abtr \abtl_{1,2} \P_{1,2}^{-1} : \A^e \otimes \M \to \M, \quad \quad  \atl = \abtr (\s \otimes \id) \P \abtl_{0,1} : \M \otimes A^e  \to \M$$
 \end{defn}

    It is straightforward to check that:

 \begin{lemma}
  \label{lemma_left_action_A^e_on_M}
  $\atr$ is a left action, and $\atl$ is a right action of $\A^e$ on $\M$, i.e. 
  $$\btr ( \id_{\A^e} \otimes \btr ) = \, \btr  ( \mu_{\A^e} \otimes \id_\M), \quad \quad \btl ( \btl \otimes \id_{\A^e}) = \btl ( \id_\M \otimes {\mu}_{\A^e})$$
  \end{lemma}
  \begin{pf} We check that 
  \begin{eqnarray}
  &\btl (\btl \otimes \id_{\A^e} ) &= \tr (\s \otimes \id) \P \tl_{0,1} \tr_{0,1} (\s \otimes \id^{\otimes 3}) \P_{0,1} \tl_{0,1}
  \nonumber\\
  &&= \tr (\s \otimes \id) \P \tr_{0,1} (\s \otimes \id^{\otimes 2}) \P_{0,1} \tl_{0,1} {\mu}_{1,2} \P_{2,3}^{-1}\nonumber\\
  &&= \tr \P {\mu}_{1,2} (\id \otimes \s \otimes \s) \P_{1,2} \tl_{0,1} {\mu}_{1,2} \P_{2,3}^{-1}\nonumber\\
  &&= \tr (\s \otimes \id) \P {\mu}_{1,2} \P_{1,2}^{-1} \tl_{0,1} {\mu}_{1,2} \P_{2,3}^{-1} = \btl ( \id_M \otimes {\mu}_{\A^e})\nonumber
  \end{eqnarray}
  where we used the ribbon property in the form ${\mu} \P (\s \otimes \s) = \s {\mu} \P^{-1}$.
  \end{pf}

So if $\A$ is a ribbon algebra, we can make any $\A$-bimodule $\M$ into both a left and a right $\A^e$-module.
 Baez performed these constructions under the assumption that $\A$ is weakly $\P$-commutative, meaning 
 ${\mu}  \P^2 = {\mu}$. 
 The ribbon automorphism $\s$ is the missing ingredient needed to make Baez's constructions work in full generality (the weakly $\P$-commutative case corresponds to $\s = \id$).
 Conversely, it is easy to check that:

\begin{lemma}\label{left_A^e_module}  If $\M$ is a left $\A^e$-module via $\btr$, then $\M$ is an $\A$-bimodule via
$$\tr : = \btr ( \id_\A \otimes \eta \otimes \id_\M) : \A \otimes \M \to \M, \quad
\tl : = \btr ( \eta \otimes \id_{\A \otimes \M})(\s \otimes \id_\M) \P : \M \otimes \A \to \M$$
If $\M$ is a right $\A^e$-module via $\btl$, then $\M$ is an $\A$-bimodule via
$$\tr : = \btl ( \id_\M \otimes \eta \otimes \id_\A) ( \id_\M \otimes \s^{-1} ) \P^{-1} : \A \otimes \M \to \M, 
\quad
\tl := \btl ( \id_{\M \otimes \A} \otimes \eta) : \M \otimes \A \to \M$$
\end{lemma}

As for any simplicial object, $\{ C_\ast (\A,\A) \}$ gives rise to a resolution of $\A$ in $\C$, called the bar resolution, when considered with the differential $\bp_{n+1} : \A^{\otimes (n+1)} \to \A^{\otimes n}$ 
 defined by $\bp_{n+1} = \sum_{i=0}^{n-1} (-1)^i {\mu}_{i,i+1}$ 
 (again, in $\A^{\otimes (n+1)}$ we number the strands 0,1, $\ldots ,n$, and ${\mu}_{i,i+1} = \id^{\otimes i} \otimes {\mu} \otimes \id^{\otimes(n-i-1)}$). 
 In general, projectivity (in $\C$) of this resolution is a subtle question. 
 However, if we take $\C = \C(\H)$ for a coquasitriangular Hopf algebra $\H$ (over $\field$) as in Section \ref{section_cqt}, then we have a forgetful functor $\C \to \kVec$, and we can 
 consider $\A^e$ simply as a $\field$-algebra.  Since $\field$ is a field $\A^e$ is projective  as a $\field$-module, and the bar resolution will then be a projective resolution in the category of modules over the $\field$-algebra $\A^e$ (in the usual sense of ring theory).
   Tensoring the bar resolution over $\A^e$ on the left by $\M$ (with right $\A^e$-module structure given by $\btl$) gives, up to isomorphism of complexes of $\field$-vector spaces, the complex defining braided Hochschild homology. Hence we have:

 \begin{thm}
\label{braided_Hoch_is_Tor}
 For $\C = \C(\H)$, with $\H$ a coquasitriangular Hopf algebra, there is an isomorphism of $\field$-vector spaces 
 $H^{\P,\s}_\ast (\A,\M) \cong \Tor^{\A^e}_\ast (\M,\A)$.\\\\
\end{thm}

\begin{figure}
\label{fig:equiv}
\[ \epsfbox{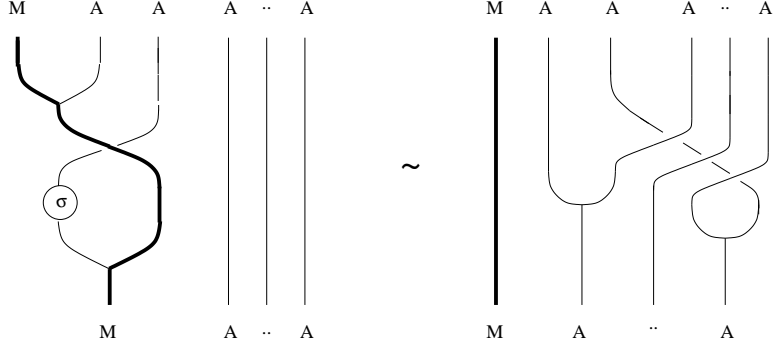} \]
\caption{The equivalence relation $\sim$  for which $\M \otimes_{\A^e} \A^{\otimes n} = \M \otimes \A^{\otimes (n+2)} / \sim$}
\label{equiv}
\end{figure}

\begin{figure}
\label{fig:HHzero}
\[ \epsfbox{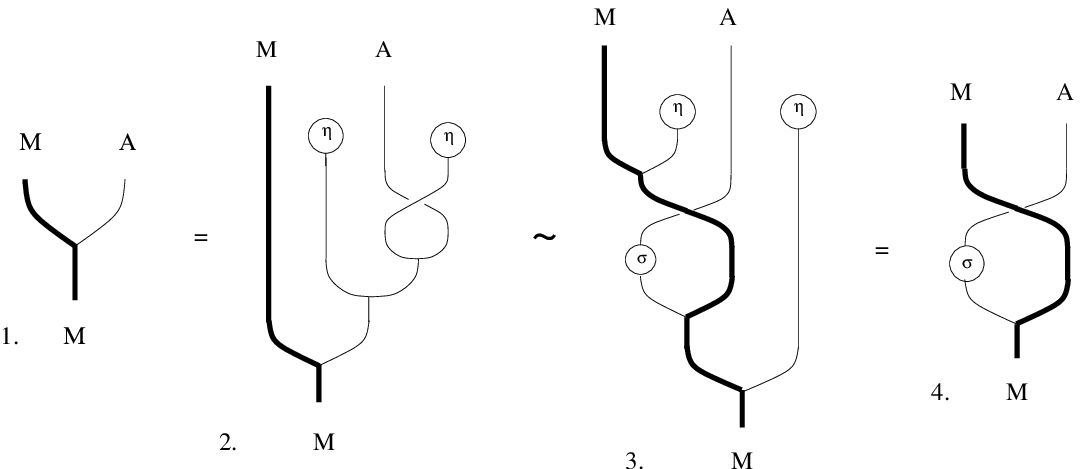} \]
\caption{Graphical proof of Theorem \ref{braided_Hoch_is_Tor} in degree zero.}
\label{HHzero}
\end{figure}

We give a graphical illustration. Figure \ref{equiv} shows the relation $\sim$ for which 
$\M \otimes_{\A^e}  \A^{\otimes n} \cong \M \otimes \A^{\otimes (n+2)} / \sim$. 
Figure \ref{HHzero} shows that $\Tor_0^{\A^e} (\M , \A) = \M \otimes_{\A^e} \A \cong \M / \{\tl \sim \tr (\s \otimes \id) \P  \}$, illustrating Theorem \ref{braided_Hoch_is_Tor} in degree zero.

  \begin{figure}
\label{fig:tn}
 \[ \epsfbox{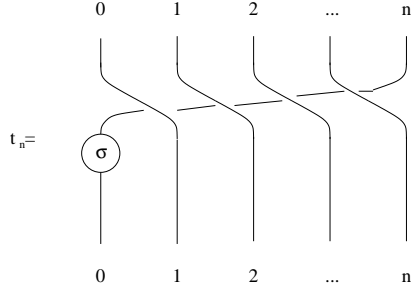} \]
\caption{The cyclic operator $t_n$}
\label{tn}
\end{figure}

\subsection{Braided cyclic homology} 
\label{section:cyclic_module}
\label{section:braidedcyclic}

Suppose now that $\M = \A$, and for each $n \geq 0$ define, as in Figure \ref{tn}, 
$$t_n = (\s \otimes \id^{\otimes n}) \P_{[0,n-1],n} : \A^{\otimes (n+1)} \to \A^{\otimes (n+1)}$$
This makes $\{ C_\ast (\A, \A) \}$ into a paracyclic object (see \cite{getzlerjones} for this notion).
Passing to the coinvariants of $T_n := t_n^{n+1}$, that is considering  the cokernels $C_n^{\P,\s} := 
\coker( \id - T_n )$ 
 we obtain the cyclic object of Akrami and Majid \cite{am}.

\begin{defn} Braided cyclic homology $HC_\ast^{\P,\s} (\A)$ is defined to be the cyclic homology of this cyclic object.
\end{defn}

\begin{example} In the setting of Example \ref{reduces_to_kmt},
 $HH^{\P,\s}_\ast (\A)$ and $HC_\ast^{\P,\s} (\A)$ reduce to Kustermans-Murphy-Tuset twisted Hochschild and cyclic homology $HH^\s_\ast (\A)$, $HC^\s_\ast (\A)$ \cite{kmt}.
 \end{example}

Note that the simplicial homology of this cyclic object, which we will denote by  $HH_\ast^{\P,\s} (\A)$, is in general different  from $H_\ast^{\P,\s} (\A,\A)$ as defined in the previous section.
The short exact sequence of complexes 
\begin{equation}
0 \to C^1_\ast \to C_\ast \to C^{\P,\s}_\ast \to 0,
 \quad C^1_\ast : = \ker \; \coker( \id - T)
\end{equation}
yields a long exact sequence
$$ \ldots \to H_n ( C^1_\ast ) \to H_n (C_\ast ) \to H_n ( C^{\P,\s}_\ast ) = HH_n^{\P,\s} (\A) \to H_{n-1} ( C^1_\ast ) \to \ldots$$
 but in general $C^1_\ast$ is not exact, as we see from the following simple example.
 
 \begin{example} Let $\A$ be the 3-dimensional unital $\bC$-algebra generated by $x$, $y$ with relations $x^2 = 0 = y^2 = xy = yx$. 
 Let $\P$ be the flip, and $\s$ the automorphism defined by $\s (x) = x$, $\s(y) = x+y$. 
 Then  $[ 1 \otimes x] \in H_1 ( C^1_\ast )$ is nontrivial.  
 \end{example}

However, in good cases the morphism $H_\ast^{\P,\s} (\A,\A) \to HH_\ast^{\P,\s} (\A)$ is an isomorphism:

 \begin{prop}\label{am_isom_baez} If $C_n \cong C_n^0
  \oplus C_n^1$ for all $n$, for $C_n^0 = \ker( \id -
  T_n)$ and $C_n^1 = \ker \; \coker (\id - T_n)$, then $HH_\ast^{\P,\s} (\A) \cong H_\ast^{\P,\s} (\A,\A)$.
 \end{prop}

See  Proposition 2.1 of \cite{hk} for a proof. 

It was shown in \cite{am} Theorem 9 
  that ribbon algebras $\A$, $\A^\f$ obtained from one another by gauge transformation (as discussed in Section \ref{cochaintwist})
   possess isomorphic cyclic (and in fact paracyclic) objects.
 Thus braided cyclic and Hochschild homology are invariant under gauge transformations.  
 
\begin{example} The standard quantised coordinate ring $\bC_q [G]$ of a complex simple Lie group $G$ (Example \ref{cqg})
 is obtained by twisting the classical coordinate ring $\bC [G]$ by the Drinfeld twist $\f$. 
  This is not a 2-cocycle, but $\partial \f$ is cocentral.
  Hence although $\bC_q [ G]$ is associative, twisting $\bC [G]$-comodule algebras produces in general nonassociative algebras \cite{am,bm}.
  
  Theorem 9 of \cite{am} thus allows us to express braided Hochschild homology of $\bC_q [G]$ as braided Hochschild homology of  $\bC [G]$. 
  Indeed, for any quasi-Hopf algebra $\H$, and any 2-cochain $\f$, we have $(\H^\cop)^{\bar{\f}} \cong (\H^\f )^\cop$. 
 Hence $(\H^\cop \otimes \H)^{\bar{\f} \otimes \f} \cong (\H^\f)^\cop \otimes \H^\f$. 
  Twisting the $\H^\cop \otimes \H$-comodule algebra $\A : = \H$ by $\bar{\f} \otimes \f$ gives 
   $$a \bullet b = a_\two b_\two \bar{\bar{\f}} ( a_\one , b_\one ) \bar\f(a_\three, b_\three) = \f (a_\one , b_\one) a_\two b_\two \bar{\f} ( a_\three , b_\three)$$
  i.e. as an algebra, 
    $\A^{\bar{\f} \otimes \f}$ is isomorphic to the twisted quasi-Hopf algebra  $\H^\f$. 
    Apply this with $\H = \bC_q [G]$, $\r$ the standard universal r-form, and $\f$ the inverse of the Drinfeld twist. 
   This gives  the braiding $\P$ as in (\ref{baez_braiding}) on $\bC_q [G]$ considered by Baez.
   Then $\H^\f \cong \bC [G]$, but the braiding obtained from $\P$ is {\bf not} the flip. 
   In particular, it does not immediately follow from \cite{am}, Theorem 9 that braided Hochschild homology of $\bC_q [G]$ can be identified with standard Hochschild homology of $\bC [G]$.
   \end{example}

\section{Braided Hochschild homology of Hopf algebras in braided categories}
\label{section:HH_Hopf}

In this Section we extend a result of Feng and Tsygan \cite{FT}, giving a simpler description of braided Hochschild homology as a derived functor in the case when $\A$ is a Hopf algebra.
 Throughout we keep the assumptions and notations of the previous Section, in particular $(\A,\s)$ will always denote a ribbon algebra in a braided monoidal category $\C$.

\begin{figure}
\label{fig:FTaction}
\[ \epsfbox{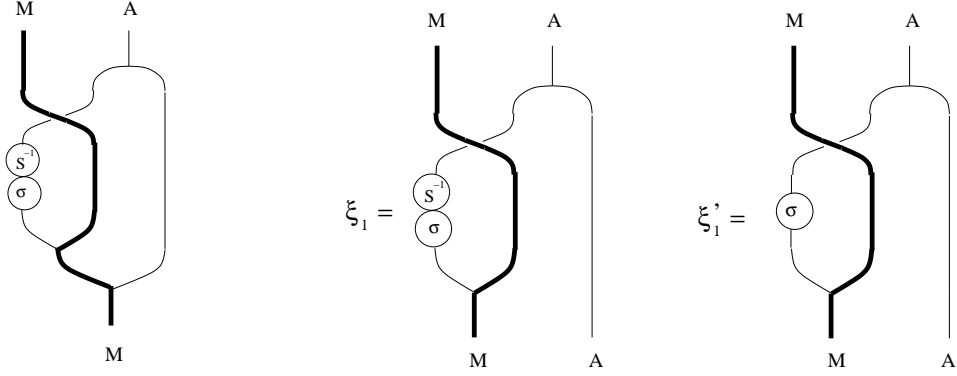} \]
\caption{The right $\A$-module structure of $\RM$, and the maps $\xi_1$, $\xip_1$}
\label{FTaction}
\end{figure}

\begin{lemma}
\label{lemma_M_is_right_A_module} Assume that $\A$ is a Hopf algebra in $\C$
 with invertible antipode  and let $\M \in \Ob(\C)$ be an $\A$-bimodule. 
Then the map
\begin{equation}
\label{black_right_action_A_on_M}
\btl \;= \tl \tr_{0,1} [ \s S^{-1} \otimes \id^{\otimes 2} ] \P_{0,1} ( \id \otimes \D)
\end{equation}
gives a right action of $\A$ on $\M$ (Figure \ref{FTaction}). We write $\RM$ for the corresponding right $\A$-module.
\end{lemma}
\begin{pf}  Writing $\D_m$ for $\id^{\otimes m} \otimes \D \otimes \id^{\otimes n}$ ($m \geq 0$), we check that 
\begin{eqnarray}
&\btl (\btl \otimes \id) &= \tl \tr_{0,1} ( \s S^{-1} \otimes \id^{\otimes 2}) \P_{0,1} \D_1 \tl_{0,1} \tr_{0,1} (\s S^{-1} \otimes \id^{\otimes 3}) \P_{0,1} \D_1
\nonumber\\
&&= \tl \tl_{0,1} \tr_{0,1} (\s S^{-1} \otimes \id^{\otimes 2} ) \P_{0,1} \P_{1,2} \tr_{0,1} ( \s S^{-1} \otimes \id^{\otimes 4}) \P_{0,1} \D_1 \D_2
\nonumber\\
&&= \tl \tr_{0,1} [ {\mu}  \P ( \s S^{-1} \otimes \s S^{-1}) \otimes \id^{\otimes 2}]  \P_{0,[1,2]} {\mu}_{3,4} \P_{2,3} \D_1 \D_2 
\nonumber\\
&&= \tl \tr_{0,1} ( \s S^{-1} \otimes \id^{\otimes 2} ) {\mu}_{0,1} \P_{0,[1,2]} {\mu}_{3,4} \P_{2,3} \D_1 \D_2
\nonumber\\
&&= \tl \tr_{0,1} ( \s S^{-1} \otimes \id^{\otimes 2}) \P_{0,1} {\mu}_{1,2} \; {\mu}_{3,4} \P_{2,3} \D_1 \D_2
\nonumber\\
&&= \tl \tr_{0,1} ( \s S^{-1} \otimes \id^{\otimes 2}) \P_{0,1} (\id \otimes \D)(\id \otimes {\mu})
= \btl (\id \otimes {\mu})
\nonumber
\end{eqnarray}
where we used the identity ${\mu} \P ( \s S^{-1} \otimes \s S^{-1} ) = \s S^{-1}  {\mu}$.
\end{pf}

Working with the coopposite Hopf algebra yields the straightforward generalisation of the action considered in \cite{FT}. 
 However the above variant is more convenient for graphical calculations.

As in Theorem \ref{braided_Hoch_is_Tor}, we now restrict to $\C = \C(\H)$, for $\H$ a coquasitriangular Hopf algebra. We note that for $\A$ a Hopf algebra in $\C$, then $\bigone_{\C(\H)} = \field$ is a left $\A$-module via the counit $\e$. 
 We now have the following generalisation of \cite{FT}, Corollary 2.5:

\begin{thm} 
\label{braided_ft} 
Let $\A$ be a Hopf algebra in  $\C(\H)$ with ribbon automorphism $\s$ and invertible antipode, and $\M$ an $\A$-bimodule. Then there is a natural isomorphism of vector spaces
$$H_n^{\P, \s} (\A,\M) \simeq \Tor_n^\A (\RM,\field)$$
where $\RM$ is the right $\A$-module of Lemma \ref{lemma_M_is_right_A_module}.
\end{thm}
\begin{pf}
 Define maps $\xi_n, \; \xip_n : \M \otimes \A^{\otimes n} \to \M \otimes \A^{\otimes n}$, $n = 1,2 , \ldots$ by
\begin{eqnarray}
&\xi_1 &= \tr_{0,1} [ \s S^{-1} \otimes \id^{\otimes 2} ] \P_{0,1} ( \id \otimes \D) \nonumber\\
&\xi_{n+1} &= \P_{[1,n],n+1}^{-1} ( \xi_1 \otimes \id^{\otimes n} )  \P_{[1,n],n+1} (\xi_n \otimes \id) \nonumber\\
&\xip_1 &= \tr_{0,1} [ \s  \otimes \id^{\otimes 2} ] \P_{0,1} ( \id \otimes \D) \nonumber\\
\label{defn_xi_n}
&\xip_{n+1} &= (\xip_n \otimes \id ) \P_{[1,n],n+1}^{-1} ( \xip_1 \otimes \id^{\otimes n} )  \P_{[1,n],n+1}
\end{eqnarray}
$\xi_1$ and $\xip_1$ are shown in Figure \ref{FTaction}. We have:

\begin{lemma} 
$\xi_n \circ \xip_n = \id = \xip_n \circ \xi_n$
\end{lemma}
\begin{pf}
We prove this  by  induction. First of all
\begin{eqnarray}
&\xip_1 \circ \xi_1 &=
\tr_{0,1} {\mu}_{0,1} (\s \otimes \s S^{-1} \otimes \id^{\otimes 2}) \P_{1,2} \P_{0,1} \P_{1,2} \D_2 \D_1
\nonumber\\
&&=
\tr_{0,1} (\s \otimes \id^{\otimes 2} ) \P_{0,1} {\mu}_{1,2} \P_{1,2}^{-2} ( \id^{\otimes 2} \otimes S^{-1} \otimes \id) \P_{1,2} \D_2 \D_1
\nonumber\\
&&=
\tr_{0,1} (\s \otimes \id^{\otimes 2} ) \P_{0,1} [ \id \otimes ( {\mu} \circ \P^{-1} ( S^{-1} \otimes \id ) \D ) \otimes \id ] \D_1
\nonumber\\
&&= \tr_{0,1} \P_{0,1} ( \id \otimes \eta \circ \e \otimes \id ) \D_1 = \id_{\M \otimes \A}
\nonumber
\end{eqnarray}
where we used the identity $\mu \circ \P^{-1} (S^{-1} \otimes \id) \D = \eta \circ \e$. 
 Suppose that $\xip_n \circ \xi_n = \id$. 
From (\ref{defn_xi_n}), 
$$\xip_{n+1} \circ \xi_{n+1} =
(\xip_n \otimes \id) \P_{[1,n],n+1}^{-1} ( \xip_1 \circ \xi_1 \otimes \id^{\otimes n} )  \P_{[1,n],n+1} (\xi_n \otimes \id ) = \xip_n \circ \xi_n \otimes \id = \id$$
The proof that $\xi_n \circ \xip_n = \id$ is completely analogous.
\end{pf}

Recall (\ref{d_j}) the maps $d_j : \M \otimes \A^{\otimes (n+1)} \to \M \otimes \A^{\otimes n}$, $0 \leq j \leq n+1$.
  We also define\\  
 ${\tilde{d}}_0,\; {\tilde{d}}_{n+1} : \M \otimes \A^{\otimes (n+1)} \to \M \otimes \A^{\otimes n}$
  by
 ${\tilde{d}}_0 = \btl \otimes \id^{\otimes n}$, ${\tilde{d}}_{n+1} = \id^{\otimes (n+1)} \otimes \e$.

 \begin{lemma} $\xip_n \circ d_0 \circ \xi_{n+1} = {\tilde{d}}_0 = \btl \otimes \id^{\otimes n}$, for all $n \geq 0$.
 \end{lemma}
 \begin{pf} First, 
 $d_0 \circ \xi_1 =  \tl \tr_{0,1} [ \s S^{-1} \otimes \id^{\otimes 2} ] \P_{0,1} ( \id \otimes \D) = \btl$, 
 by  definition of $\btl$ (\ref{black_right_action_A_on_M}).
  Next, 
  \begin{eqnarray}
  & \xip_1 \circ d_0 \circ \xi_2
   &= \xip_1\;  \tl_{0,1} \P_{1,2}^{-1} \; ( \xi_1 \otimes \id ) \P_{1,2} ( \xi_1 \otimes \id)
    \nonumber\\
   && = 
 \tr_{0,1} \; ( \s \otimes \id^{\otimes 2} ) \P_{0,1} \D_1 \tl_{0,1} \; \P_{1,2}^{-1} \tr_{0,1} ( \s S^{-1} \otimes \id^{\otimes 3} ) \P_{0,1} \D_1
  \P_{1,2} \tr_{0,1} ( \s S^{-1} \otimes \id^{\otimes 3}) \P_{0,1} \D_1  
   \nonumber\\
 && = 
 \tl_{0,1} \P_{1,2}^{-1} \tr_{0,1} (\s \otimes \id^{\otimes 3}) \P_{0,1} \D_1 \tr_{0,1} ( \s S^{-1} \otimes \id^{\otimes 3}) \P_{0,1} \D_1 
  \tr_{0,1} ( \s S^{-1} \otimes \id^{\otimes 3}) \P_{0,1} \P_{2,3} \D_1
 \nonumber\\
 && = \tl_{0,1} \P_{1,2}^{-1} \tr_{0,1} ( \s S^{-1} \otimes \id^{\otimes 3}) \P_{0,1} \P_{2,3} \D_1
   = \tl_{0,1} \tr_{0,1} (\s S^{-1} \otimes \id^{\otimes 3} ) \P_{0,1} \D_1 = \btl \otimes \id
 \nonumber
  \end{eqnarray}
 By induction, if $\xip_n \circ d_0 \circ \xi_{n+1} = \btl \otimes \id^{\otimes n}$, then 
  \begin{eqnarray}
  &\xip_{n+1} \circ d_0 \circ \xi_{n+2} &= 
  (\xip_n \otimes \id) \P_{[1,n],n+1}^{-1} (\xip_1 \otimes \id^{\otimes n}) \P_{[1,n], n+1} \tl_{0,1}  
  \P_{[1,n], n+1}^{-1} (\xi_1 \otimes \id^{\otimes n}) \P_{[1,n+1],n+2} ( \xi_{n+1} \otimes \id) \nonumber\\
  && =  (\xip_n \otimes \id) \P_{[1,n],n+1}^{-1} (\xip_1 \otimes \id^{\otimes n}) \P_{[1,n],n+1} \tl_{0,1} \P_{[1,n+1], n+2}^{-1}  \tr_{0,1} ( \s S^{-1} \otimes \id^{\otimes n+2} ) \P_{0,1}\nonumber\\
  &&\quad \quad  \D_1 \P_{[1,n+1], n+2} ( \xi_{n+1} \otimes \id)\nonumber\\ 
  &&= ( \xip_n \otimes \id) \P_{[1,n],n+1}^{-1} (\xip_1 \otimes \id^{\otimes n}) (\xi_1 \otimes \id^{\otimes n}) \P_{[1,n],n+1} \tl_{0,1} (\xi_{n+1} \otimes \id)\nonumber\\
  && = ( \xip_n \otimes \id) \tl_{0,1} ( \xi_{n+1} \otimes \id) = \btl \otimes \id^{\otimes n+1}\nonumber
  \end{eqnarray}
 \end{pf}

  \begin{lemma} $\xip_n \circ d_1 \circ \xi_{n+1} = d_1$ for all $n \geq 1$.
\end{lemma}
\begin{pf} 
First of all, 
\begin{eqnarray}
& \xip_1 \circ d_1 \circ \xi_2 &= \tr_{0,1} ( \s \otimes \id^{\otimes 2}) ( \id \otimes \D \circ {\mu}) \P_{1,2}^{-1} \tr_{0,1} ( \s S^{-1} \otimes \id^{\otimes 3}) \P_{0,1} \D_1 \P_{1,2} (\xi_1 \otimes \id)
\nonumber\\
&& = \tr_{0,1} (\s \otimes \id^{\otimes 2} ) \P_{0,1} {\mu}_{1,2} {\mu}_{3,4} \P_{2,3} \D_1 \D_2 \P_{1,2}^{-1} \tr_{0,1} ( \s S^{-1} \otimes \id^{\otimes 3}) \P_{0,1} \D_1 \P_{1,2} (\xi_1 \otimes \id)\nonumber\\
&& = \tr_{0,1} \P_{0,1} {\mu}_{1,2} (\id \otimes \s \otimes \s S^{-1} \otimes \id) \P_{1,2} {\mu}_{1,2} {\mu}_{3,4} \P_{3,4} \D_2 \D_3 \P_{1,2} \D_2 (\xi_1 \otimes \id)\nonumber\\
&& = \tr_{0,1} (\s \otimes \id^{\otimes 2}) \P_{0,1} {\mu}_{1,2} \P_{1,2}^{-1} (\id \otimes S^{-1} \otimes \id^{\otimes 2}) {\mu}_{2,3} {\mu}_{4,5} \P_{3,4} \D_2 \D_3 \P_{1,2} \D_2 (\xi_1 \otimes \id)\nonumber\\
&& = \tr_{0,1} (\s \otimes \id^{\otimes 2} ) \P_{0,1} {\mu}_{1,2} \P_{2,3} {\mu}_{1,2} {\mu}_{3,4} \P_{2,3} (\id^{\otimes 5} \otimes S^{-1} ) \D_1 \D_2 \P_{2,3}^{-1} \D_2 (\xi_1 \otimes \id)\nonumber\\
&& = \tr_{0,1} (\s \otimes \id^{\otimes 2} ) \P_{0,1} {\mu}_{1,2} {\mu}_{2,3} \P_{3,4} {\mu}_{3,4} \P_{2,3} (\id^{\otimes 5} \otimes S^{-1} ) \D_1 \D_2 \P_{2,3}^{-1} \D_2 (\xi_1 \otimes \id)\nonumber\\ 
&& = \tr_{0,1} (\s \otimes \id^{\otimes 2} ) \P_{0,1} {\mu}_{1,2} {\mu}_{2,3} \P_{3,4} {\mu}_{3,4} \P_{2,3} \P_{4,5}^{-1} \P_{3,4}^{-1} (\id^{\otimes 3} \otimes S^{-1} \otimes \id^{\otimes 2}) \D_1 \D_2  \D_2 (\xi_1 \otimes \id)\nonumber\\
&& = \tr_{0,1} (\s \otimes \id^{\otimes 2} ) \P_{0,1} {\mu}_{1,2}  {\mu}_{3,4} \P_{2,3} \D_1 ( \id^{\otimes 2} \otimes [ {\mu} \P^{-1} ( S^{-1} \otimes \id) \D ] \otimes \id) \D_2 (\xi_1 \otimes \id)\nonumber\\
&& = \tr_{0,1} (\s \otimes \id^{\otimes 2} ) \P_{0,1} {\mu}_{1,2}  {\mu}_{3,4} \P_{2,3} \D_1 (\id^{\otimes 2} \otimes \eta \circ \e \otimes \id) \D_2 (\xi_1 \otimes \id)\nonumber\\
&& = {\mu}_{1,2} \tr_{0,1} (\s \otimes \id^{\otimes 3}) \P_{0,1} \D_1 (\xi_1 \otimes \id)
= {\mu}_{1,2} ( \xip_1 \otimes \id)(\xi_1 \otimes \id) =  {\mu}_{1,2} = d_1\nonumber
\end{eqnarray}

Now suppose that $\xip_n \circ d_1 \circ \xi_{n+1} =d_1$. Then
\begin{eqnarray}
& \xip_{n+1} \circ d_1 \circ \xi_{n+2} &= 
(\xip_n \otimes \id ) \P_{[1,n],n+1}^{-1} ( \xip_1 \otimes \id^{\otimes n} )  \P_{[1,n],n+1} d_1 
 \P_{[1,n+1],n+2}^{-1} ( \xi_1 \otimes \id^{\otimes n+1} )\nonumber\\
 &&\quad \quad  \P_{[1,n+1],n+2} (\xi_{n+1} \otimes \id) \nonumber\\
 &&= (\xip_n \otimes \id ) \P_{[1,n],n+1}^{-1} ( \xip_1 \otimes \id^{\otimes n} ) (\xi_1 \otimes \id^{\otimes n}) {\mu}_{2,3} \P_{[1,n+1], n+2} ( \xi_{n+1} \otimes \id)\nonumber\\
  &&= (\xip_n \otimes \id ) \P_{[1,n],n+1}^{-1} {\mu}_{2,3} \P_{[1,n+1], n+2} ( \xi_{n+1} \otimes \id)\nonumber\\ &&= (\xip_n \otimes \id ) {\mu}_{1,2} ( \xi_{n+1} \otimes \id)
    = (\xip_n \circ d_1 \circ \xi_{n+1}) \otimes \id = d_1
    \nonumber
\end{eqnarray}
\end{pf}

\begin{lemma}
$\xip_n \circ d_{n+1} \circ \xi_{n+1} = ( \id^{\otimes n+1} \otimes \e ) = {\tilde{d}}_{n+1}$ for all $n \geq 1$.
\end{lemma}
\begin{pf} We have  
\begin{eqnarray}
& d_1 \circ \xi_1 & = \tr_{0,1} (\s \otimes \id) \P_{0,1}  \tr_{0,1} ( \s S^{-1} \otimes \id^{\otimes 2} ) \P_{0,1} \D_1 \nonumber\\
&& = \tr_{0,1} (\s \otimes \id) \P_{0,1} {\mu}_{1,2} \P_{1,2} (\id \otimes S^{-1} \otimes \id) \D_1 = \id \otimes \e\nonumber\\
&\Rightarrow \; \xip_{n} \circ d_{n+1} \circ \xi_{n+1} & =  \xip_{n} \tr_{0,1} ( \s \otimes \id^{\otimes n+1}) \P_{[0,n],n+1}  \P_{[1,n],n+1}^{-1} ( \xi_1 \otimes \id^{\otimes n} )  \P_{[1,n],n+1} (\xi_{n} \otimes \id) \nonumber\\
&&= \xip_{n} ( d_1 \circ \xi_1 \otimes \id^{\otimes n} ) \P_{[1,n],n+1} ( \xi_{n} \otimes \id) 
\nonumber\\
&& = \xip_{n} ( \id \otimes \e \otimes \id^{\otimes n}) \P_{[1,n], n+1} ( \xi_{n} \otimes \id) = \id^{\otimes  n+1} \otimes \e\nonumber
\end{eqnarray}
for all $n \geq 1$. 
\end{pf}

Finally it is straightforward to check that:

\begin{lemma} $\xip_n \circ d_i \circ \xi_{n+1} = d_i$ for $1 \leq i \leq n$.
\end{lemma}

Therefore $\xi_n$, $\xip_n$ define isomorphisms of complexes between the braided Hochschild complex and the complex $\{ \M \otimes \A^{\otimes n}, {\tilde{b}}_n \}$, where ${\tilde{b}}_{n+1} = \tilde{d_0} + \sum_{j=1}^{n} (-1)^j d_j + {\tilde{d}}_{n+1} : M \otimes \A^{\otimes (n+1)} \to M \otimes \A^{\otimes n}$. 
That the homology of the latter complex is $\Tor^\A_\ast ( R(M), \field)$ follows from \cite{ce} Cor. IX.4.4, 
applied with $K$, $\Lambda$, $\A$, $\B$ being $\field$, $\A$, $\field$, $R(M)$ respectively.
\end{pf}

\begin{figure}
\label{fig:FTHHzero}
\[ \epsfbox{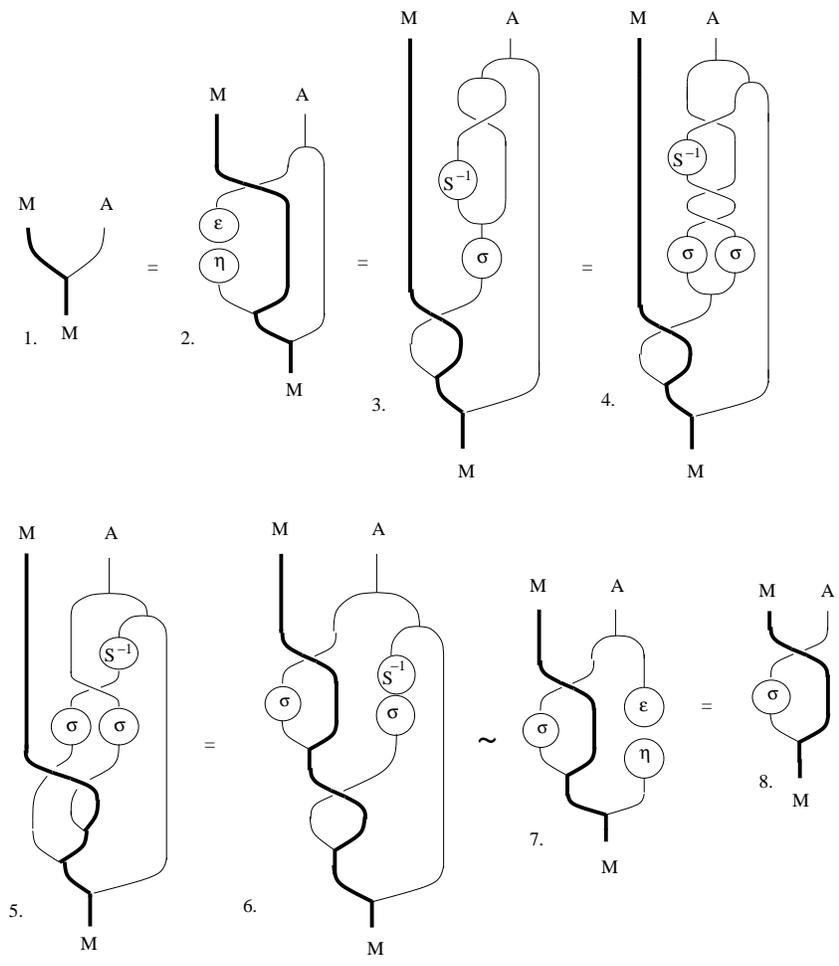} \]
\caption{Graphical proof of Theorem \ref{braided_ft} in degree zero.} 
\label{FT_HHzero}
\end{figure}

 Figure \ref{FT_HHzero} shows a graphical proof in degree zero.

\section{The braided line and the braided plane}
\label{section:braided_line}

In this section we consider the  braided line and braided plane. 
 We observe that braided Hochschild homology is  less degenerate than ordinary Hochschild homology of the classical counterparts.
  Throughout this section  $q \in \bC$ will denote a nonzero parameter which is not a root of unity.

\subsection{The braided line}
\label{braided_line}

 Let $H$ be the commutative Hopf algebra $\field \bZ = \field [t,t^{-1}]$, with coproduct $\D(t) = t \otimes t$. 
 Then $H$ is coquasitriangular, via $\r( t^m , t^n ) = q^{mn}$ for all $m$, $n \in \bZ$. 
  Let $\field [x]$ be the unital $\field$-algebra in a single indeterminate $x$. 
  Then $\field [x]$ is a right $H$-comodule algebra, via $x^m \mapsto x^m \otimes t^m$, for all $m \geq 1$. 
 The braiding  (\ref{rightcomodbraiding}) is
$\P ( x^m \otimes x^n ) = q^{mn} x^n \otimes x^m$. 
The braided line $\A$ is  the braided Hopf algebra given by $\field [x]$ with    
$$\D(1) = 1 \otimes 1, \quad \D(x) = x \otimes 1 + 1 \otimes x, \quad \varepsilon (x) = 0, \quad S(x) = -x$$
It follows from the ribbon relation $\mu \circ ( \s \otimes \s) \circ \Psi^2 = \s \circ {\mu}$
that any ribbon automorphism satisfies  $\s(1) =1$, $\s(x^n) = q^{n(n-1)} \s(x)^n$ for all $n \geq 1$.
 Taking $\s(x) = \sum_{n \geq 0} \l_n x^n$,  compatibility with the braiding forces $\l_n =0$ for $n \neq 1$, so $\s(x) = \l x$, $\s(x^n) = ( \l q^{n-1})^n x^n$ for all $n \geq 1$.
 Since we need $\s$ to be invertible,  $\l \neq 0$.

  Every elementary tensor $x^{m_0} \otimes \ldots \otimes x^{m_n}$ is an eigenvector of $T_n$, hence Proposition \ref{am_isom_baez} applies and $HH_n^{\P,\s} (\A) \cong H_n^{\P,\s} (\A,\A)$ for all $n$.

\begin{prop} We have $HH^{\Psi, \s}_n ( \A) = 0$ for all $n \geq 2$, for all $q$, $\l$.
\begin{enumerate}
\item For $\l \notin q^{- \bN}$, then $HH^{\Psi, \s}_0 ( \A) = \field [1] $, $HH^{\Psi, \s}_1 ( \A) =0$.
\item For $\l = q^{- N}$, some $N \in \bN$,  then  $HH^{\Psi, \s}_0 ( \A) = \field [1] \oplus \field [ x^{N+1}]$, $HH^{\Psi, \s}_1 ( \A) = \field [ x^N \otimes x]$.
\end{enumerate}
\end{prop}
\begin{pf} We use Theorem \ref{braided_ft}. 
 Consider the resolution of $\field$ by free left $\A$-modules:
$0 \to \A \to^\varphi \A \to^\e \field \to 0$, where $\varphi(x^n ) = x^{n+1}$. 
Tensoring on the left by $R(\A) \otimes_\A -$, where $R(\A)$ has underlying space $\A$ and right action
$$x^n \btl x = x^{n+1} + q^n \s(S^{-1} (x)) x^n   = (1 - \l q^n ) x^{n+1}$$
 gives $HH^{\Psi, \s}_\ast ( \A)$ as the homology of 
$0 \to R(\A) \otimes_\A \A \to^\phi R(\A) \otimes_\A \A$. 
Here 
$$\phi( x^n \otimes 1) = x^n \otimes x = x^n \btl x \otimes 1 = (1 - \l q^n ) x^{n+1} \otimes 1$$
 It follows that if $\l \notin q^{-\bN}$, then $HH^{\Psi, \s}_0 ( \A) = \field [1] $, $HH^{\Psi, \s}_1 ( \A) =0$, whereas if $\l = q^{-N}$, then $\ker(\phi)$ is one-dimensional, and $\coker( \phi)$ is two-dimensional, with generators as  above. 
\end{pf}

Using standard spectral sequence arguments \cite{hk} we calculate the braided cyclic homology. 
 We note that for the map $B_0 : HH^{\Psi, \s}_0 ( \A) \to HH^{\Psi, \s}_1 ( \A)$ we have  $B_0 [ x^{n+1} ] = (n+1) [ x^n \otimes x]$.

\begin{cor} For $q$ not a root of unity, we have $HC^{\Psi, \s}_{2n+1} ( \A) = 0$,\\ 
$HC^{\Psi, \s}_{2n+2} ( \A) = \field [1]$, for $n \geq 0$, and 
$HC^{\Psi, \s}_0 ( \A) = \left\{ 
\begin{array}{cc}
\field [1] & : \l \notin q^{-\bN} \\
\field [1] \oplus \field [ x^{N+1}] & : \l = q^{-N}
\end{array} \right.$
\end{cor}

\begin{remark} In the classical case $\l = 1 =q$,  
 $HH_n ( \A ) = 0$ for  $n \geq 2$, and 
 $HH_1 ( \A)$, $HH_0 ( \A)$ are both infinite dimensional, spanned by $\{ [ x^n \otimes x] \}_{n \geq 0}$, $\{ [x^n ] \}_{n \geq 0}$. Hence $HC_0 ( \A )$ is infinite dimensional,  
 $HC_{2n+1} ( \A ) = 0$ and $HC_{2n+2} ( \A ) = \field [ 1]$, for $n \geq 0$.
\end{remark}

\subsection{The braided plane}\label{section_braided_plane}

In this section, $\A$ denotes 
Manin's quantum plane, that is, 
the unital algebra (over $\field$) 
generated by indeterminates $x$, $y$ 
satisfying $yx = qxy$. 
This is a 
right $\H = \bC_q [SL(2)]$-comodule 
algebra via
$$[y\;, \; x ] \mapsto [y\;, \; x ] \otimes 
\left[
\begin{array}{cc}
a & b \cr
c & d \cr
\end{array}
\right]
= [\; y \otimes a + x \otimes c\;, y \otimes b + x \otimes d\; ]
$$
One can identify $\A$ with the subalgebra of $\H$
generated by $a,b$. Under this identification, 
$y$ becomes $a$, $x$ becomes $b$, and the coaction
becomes $ \Delta $ (so $\A$ is an ``embeddable quantum
space'' of $\H$). Obviously, the monomials $x^my^n$
form a vector space basis of $\A$. Under the embedding
into $\H$ these monomials become proportional to the
Peter-Weyl basis elements $C^{(m+n)}_{0m}$. The precise
factor is given e.g.~in \cite{KS}, but it is
irrelevant for us.

The braiding (\ref{rightcomodbraiding}) induced 
from the coquasitriangular structure of $\H$ is determined by: 
\begin{eqnarray}
&&\P ( x \otimes x ) = q^{1/2} x \otimes x,\quad \P ( x \otimes y) = q^{-1/2}  y \otimes x\nonumber\\
&&\P ( y \otimes x) = q^{-1/2}[ x \otimes y + (q - q^{-1}) y \otimes x],\quad \P ( y \otimes y) = q^{1/2} y \otimes y\nonumber\\
&\Rightarrow& \P( x^m y^n \otimes x) = q^{(m-n)/2} x \otimes x^m y^n + q^{(n-m-2)/2} f(n) y \otimes x^{m+1} y^{n-1},\nonumber\\
&&\P( x^m y^n \otimes y) = q^{(n-m)/2} y \otimes x^m y^n,\quad
\P( x \otimes x^m y^n ) = q^{(m-n)/2} x^m y^n \otimes x,\nonumber\\
&&\P( y \otimes x^m y^n ) = q^{(n-m)/2} x^m y^n \otimes y + q^{(m-n-2)/2} f(m) x^{m-1} y^{n+1} \otimes x\nonumber
\end{eqnarray}
 where $f(n) = q^n - q^{-n}$.  
 There is a braided Hopf algebra structure on $\A$ given by \cite{shahn_book}:
$$\D(x) = x \otimes 1 + 1 \otimes x, \quad \D( y ) = y \otimes 1 + 1 \otimes y, \quad \e(x) = 0 = \e(y), \quad S(x) = -x, \quad S(y) = -y$$
\begin{prop}
The ribbon automorphisms of $\A$ 
are parametrised by $\l \neq 0$ and given by 
$$
\s( x^m y^n ) = \l^{m+n} q^{(m+n)(m+n-1)/2}  x^m y^n.
$$ 
\end{prop}
\begin{pf} Applying the ribbon relation to  $\s(yx - qxy)$ gives 
 $\s(y) \s(x) = q \s(x) \s(y)$. Furthermore, one easily
shows by induction on $m+n$ that $ \s(x^my^n)$ 
is a linear combination of monomials in $\s(x),\s(y)$ 
(using the above formulae for the braiding). That is, 
$\s(y),\s(x)$ generate $\A$ as an algebra.
It follows that $\s(x) = \l x$, $\s(y) = \rho y$, for some $\l$, $\rho \in \bC$.
 Compatibility with the braiding forces $\l = \rho$.
Using the defining property, one easily 
checks that $ \s$ extends to $ x^m$ and 
$y^n$ by the formula given. To derive the formula in
 general one can proceed as follows: Twist the 
$\H$-comodule algebra $\A$ by the 2-cocycle 
$\q=\r_{21}\r$, that is, consider the new product
$f \bullet g:=f_{(0)}g_{(0)}\bar\q(f_{(1)},g_{(1)})$. 
The defining property of 
the ribbon automorphism 
is $ \s (f \bullet g)=\s(f) \s(g)$, and in
particular $\s(x^m \bullet y^n)=
\lambda^{m+n}
q^{m(m-1)/2+n(n-1)/2}x^my^n$.  
But since
$\q=\partial \varphi^2$ with $\varphi^2$ given in
(\ref{esistspaet}), we have also an isomorphism 
$f \mapsto \eta (f):=f_{(0)} \varphi^2(f_{(1)})$ between the two
products, and using the explicit formula for 
$ \varphi^2$ and the fact that $x^my^n$ is proportional
 to $C^{(m+n)}_{0m}$, we obtain
\begin{eqnarray}
&&\eta (x^m) \bullet \eta (y^n)=\eta (x^my^n) \nonumber\\
& \Leftrightarrow &
	(-1)^{m} q^{-(m^2/2+m)}x^m \bullet
	(-1)^{n} q^{-(n^2/2+n)}y^n=
	(-1)^{m+n} q^{-((m+n)^2/2+m+n)}x^my^n \nonumber\\ 
& \Leftrightarrow &
	x^m \bullet y^n=
	q^{-mn}x^my^n. \nonumber 
\end{eqnarray}
The claim follows.
\end{pf}

In other words, the ribbon automorphisms 
arise as $\s(f)=f_{(0)} \mathbf{s}_\l(f_{(1)})$, where 
$\mathbf{s}_\l(C^{m}_{rs})=\l^{m} q^{m(m-1)/2}
\delta_{rs}$. It is now easily shown by induction that 
$T_n$ is in fact the ribbon automorphism 
of $\A^{\otimes n+1}$, that is, is given by 
$T_n(f \otimes g \otimes \cdots \otimes h)=
f_{(0)} \otimes g_{(0)} \otimes \cdots \otimes h_{(0)}
\mathbf{s}_\l(f_{(1)}g_{(1)}\cdots h_{(1)})$. In
particular, $T_n$ acts by scalar multiplication on the
irreducible subcomodules of $\A^{\otimes n+1}$.
Hence $T_n$ is diagonalisable, and  
by Proposition \ref{am_isom_baez} we have:
  
  \begin{lemma}\label{diag} $HH_\ast^{\P,\s} (\A) = H_\ast^{\Psi,\s} (\A,\A)$.
  \end{lemma}

\begin{prop} 1. If $\l \notin q^{-\bN/2}$, then $HH_0^{\Psi,\s} ( \A )  = \bC [1]$, $HH_1^{\Psi,\s} ( \A ) = 0 = HH_2^{\Psi,\s} ( \A ) $.\\
2. If $\l = q^{-N/2}$, some $N \in \bN$, then $HH_0^{\Psi,\s} ( \A )  \cong \bC^{N+3}$, $HH_1^{\Psi,\s} ( \A ) \cong \bC^{2N+2}$,  $HH_2^{\Psi,\s} ( \A ) \cong \bC^{N+1}$
\end{prop}
\begin{pf} We compute $H^{\Psi,\s}_* (\A, \A)$ via an explicit resolution, then by Lemma \ref{diag} we  identify this with $HH^{\P,\s}_\ast ( \A)$. 
 Let $\M_2$, $\M_1$ be the free left $\A$-modules with bases 
$\{ [x \wedge y] \}$,  $\{ [x], [y] \}$ respectively, and let $\M_0 = \A$. 
Define left $\A$-module maps $f_i : \M_i \to \M_{i-1}$, $i = 2,1$ by
$$f_2  [x \wedge y]  = y [x] - q x [y], \quad f_1 [x] = x, \quad f_1 [y]= y$$
Then  the sequence 
$0 \to \M_2 \to^{f_2} \M_1 \to^{f_1} \M_0 \to^\e \bC$
is a free resolution of $\bC$ ($\A$ is a Koszul algebra, and this is the Koszul resolution of the trivial $\A$-module $\bC$). 
We tensor on the left by $R(\A)$,  the right $\A$-module with underlying space $\A$ and right action
\begin{eqnarray}
\label{right_action_x_y}
&& x^m y^n \btl x = q^n (1 - \l q^{(m+n)/2}) x^{m+1} y^n, \quad
 x^m y^n \btl y = (1 - \l q^{(m+n)/2}) x^{m} y^{n+1}
\end{eqnarray}
So  $HH^{\Psi,\s}_* (\A)$ is the homology of the complex
$0 \to \A \to^{g_2} \A \oplus \A \to^{g_1} \A \to 0$
where 
$$g_2 (a) = ( a \btl y , - q a \btl x), \quad g_1 (a,b) = a \btl x + b \btl y$$

 Using (\ref{right_action_x_y}),  then for $\l \notin q^{-\bN /2}$ we have $\ker(g_2) = 0$, $\im(g_2) = \ker(g_1)$, and $A / \im ( g_1) = \bC [1]$. 
 For $\l = q^{-N/2}$ ($N \in \bN$), then 
 $\ker(g_2) = \span \{ x^m y^n \; | \; m+n = N \; \} \cong \bC^{N+1}$, 
 $\ker(g_1) / \im(g_2) = \span \{ (x^m y^n ,0), (0, x^m y^n ) \; | \; m+n = N \; \} \cong \bC^{2N+2}$, and 
 $\A / \im(g_1) = k[1] \oplus \span \{ x^m y^n \; | \; m+n=N+1 \; \} \cong \bC^{N+3}$. 
  We can identify the generators of $HH_1^{\Psi,\s} ( \A )$ with the 1-cycles $x^m y^n \otimes x$, $x^m y^n \otimes y$, for $m+n = N$. 
   \end{pf}

\begin{thm} 1. For $\l \notin q^{- \bN /2}$, then $HC^{\P,\s}_{2n} (\A) \cong \bC$, generated by $[1]$, and $HC_{2n+1}^{\P,\s} (\A) =0$.\\
2. For $\l = q^{-N/2}$, 
$HC_0^{\P,\s} ( \A) \cong \bC^{N+3}$, $HC_1^{\P,\s} (\A)  \cong \bC^N$,
$HC_{2n+2}^{\P,\s} (\A) \cong \bC^2$, 
$HC_{2n+3}^{\P,\s} (\A) =0$, for all $n \geq 0$.
\end{thm}
\begin{pf} We calculate $HC^{\P,\s}_\ast (\A)$ as total homology of the mixed $(B,b)$-bicomplex associated to the cyclic object of Section \ref{section:braidedcyclic}, as in \cite{hk}.
 For $\l \notin q^{-\bN /2}$, the spectral sequence stabilises at the first page. 
 For $\l = q^{-N/2}$ we need to calculate with the maps $B_i : HH^{\P,\s}_i ( \A) \to HH^{\P,\s}_{i+1} (\A)$, $i=0,1$:
 $$B_0 [a] = [ 1 \otimes a], \quad B_1 [ a \otimes b] = [ 1 \otimes ( (\id - t_1 ) ( a \otimes b) )]$$
 In the same way as in \cite{hk}, Lemma 2.2, we have
 $$B_0 [ x^{N+1} ] = (N+1) [ x^N \otimes x],\quad B_0 [ x^m y^n ] = n [ x^m y^{n-1} \otimes y ] + m q^{-n} [ x^{m-1} y^n \otimes x ], \quad m+n = N+1$$
Since $b_2 ( 1 \otimes 1 \otimes 1) = 1 \otimes 1 = B_0 (1)$, we have $\ker (B_0) = \bC [1]$, $\im (B_0) \cong \bC^{N+2}$, hence $HC_1^{\P,\s} (\A) \cong \bC^N$, with generators the (equivalence classes of the) elements $x^j y^{N-j} \otimes x$, equivalently $x^{j+1} y^{N-j-1} \otimes y$, for $j = 0, 1, \ldots ,N-1$.
 Finally we show that $\im ( B_1) \cong \bC^N$. We have 
\begin{equation}\label{B1}
B_1 ( x^m y^n \otimes x) = 1 \otimes x^m y^n \otimes x - q^{-n} (1 \otimes x \otimes x^m y^n ) - q^{-m-1} f(n) (1 \otimes y \otimes x^{m+1} y^{n-1})
\end{equation}
Consider the linear functional $\tau_{s,t} : \A \to \bC$ defined by 
$\tau_{s,t} ( x^i y^j ) = \d_{s,i} \d_{t,j}$.
 Then for $s+t = N+1$, $\tau_{s,t}$ is a nontrivial braided cyclic 0-cocycle on $\A$. In particular 
 $\tau_{s,t} \circ b_1 =0$, where $b_1 : \A^{\otimes 2} \to \A$ is defined by $b_1 =  \mu - \mu \circ ( \s \otimes \id) \circ \P$. 
 Let $\partial_1$, $\partial_2 : \A \to \A$ be the derivations defined by 
 $$\partial_1 (x) =x, \quad 
\partial_1 (y) = 0, \quad \partial_2 (x) = 0, \quad \partial_2 (y) =y$$
and extended by $\partial_i (ab) = \partial_i (a) b + a \partial_i (b)$. Then 
$$\phi_{s,t} : \A^{\otimes 3} \to \bC, \quad \phi_{s,t} (a \otimes b \otimes c) = \tau_{s,t} ( \; a \; [ \; \partial_1 (b) \; \partial_2 (c) - \partial_2 (b) \; \partial_1 (c)\; ] \; )$$
is a braided Hochschild 2-cocycle, meaning $\phi_{s,t} \circ b_3 =0$. Using (\ref{B1}), 
$$\phi_{m+1,n} ( B_1 ( x^m y^n \otimes x) ) = (  (m-n+1) q^n - (m+n+1) q^{-n} )  \tau_{m+1,n} (x^{m+1} y^n)$$
which is nonzero for $n \neq 0$. Hence $B_1 ( x^m y^n \otimes x)$ represent distinct nontrivial elements of  $HH_2^{\P,\s} (\A)$, for $n =1,2, \ldots, N$. 
Since $HH_1^{\P,\s} (\A) / \im( B_0) \cong \bC^N$, it follows that $\im (B_1) \cong \bC^N$. 
 Hence $HH_2^{\P,\s} (\A) / \im(B_1) \cong \bC$, and $\ker (B_1) = \im (B_0)$. 
The spectral sequence stabilises at the second page, and we have:
\begin{eqnarray}
&&HC_0^{\P,\s} ( \A) = HH_0^{\P,\s} (\A) \cong \bC^{N+3},\quad HC_1^{\P,\s} (\A) = HH_1^{\P,\s} (\A) / \im (B_0) \cong \bC^N,\nonumber\\
&&HC_{2n+2}^{\P,\s} (\A) = \ker (B_0) \oplus HH_2^{\P,\s}(\A) / \im (B_1) \cong \bC^2, \quad
HC_{2n+3}^{\P,\s} (\A) = \ker(B_1) / \im(B_0) =0,\nonumber
\end{eqnarray}
for all $n \geq 0$. 
\end{pf}

\section{Braided $SL(2)$}
\label{qg}

\subsection{Braided Hopf algebras associated to coquasitriangular Hopf algebras}
\label{section_braidedhopf}
Let $(\H,\r)$ be a coquasitriangular Hopf algebra.
 Let $\A \in \Ob( \C(\H))$ be $\H$ equipped with the right adjoint coaction 
 $\Ad_R (a) =   a_\two \otimes S( a_\one ) a_\three$. 
 In general $\A$ is not a right $\H$-comodule algebra.
  Now define $\B = \B(\H)$ to be the algebra with underlying vector space $\H$ and (associative) product
\begin{equation}
\label{B(A)_product}
a * b =  a_\two b_\three \; \r ( a_\one , b_\two ) \r( a_\three , S b_\one ) = a_\two b_\two \; \r ( S( a_\one) a_\three , S b_\one )
\end{equation}
Then $B$ is via $\Ad_R$ a right $\H$-comodule algebra  and in fact a Hopf algebra in $\C(\H)$ with coproduct and antipode given by
\begin{equation}\label{B(A)_antipode}
\D(a) = a_\one \otimes a_\two, \quad \underline{S} (a) =  S( a_\two ) \r ( S^2 ( a_\three ) S( a_\one ) , a_\four )
\end{equation}
(see \cite{shahn_book} and  \cite{KS} Section 10.3.2). 
We call $B$ the braided Hopf algebra associated to $\H$, alternatively the transmutation of $\H$.
  The coaction $\Ad_R$  gives a braiding 
 \begin{equation}
\label{B(A)_braiding}
\P_{\B} : \B \otimes \B \to  \B \otimes \B, \quad \P_{\B} ( a \otimes b) =  b_\two \otimes a_\two \; \r( S( a_\one) a_\three , S( b_\one ) b_\three)
\end{equation}

\subsection{Braided homology of quantum $SL(2)$}
\label{section:SLq2}

Our aim in this Section is to apply this to $\A = \bC_q [ SL(2)]$  as defined in (\ref{slq2_one},\ref{slq2_two}).
 The universal r-form was explicitly recalled in (\ref{slq2_rform}).
  The resulting canonical braiding (\ref{canonical_braiding}) is defined on generators by:
\begin{eqnarray}
&\P(a \otimes a) = q^{1/2} a \otimes a,\quad& \P(a \otimes b) = q^{-1/2} b \otimes a + q^{-1/2} (q - q^{-1}) a \otimes b\nonumber\\
&\P(a \otimes c) = q^{1/2} c \otimes a,\quad& \P(a \otimes d) = q^{-1/2} d \otimes a + q^{-1/2}(q - q^{-1}) c \otimes b\nonumber\\
&\P(b \otimes a) = q^{-1/2} a \otimes b, \quad& \P(b \otimes b) = q^{1/2} b \otimes b\nonumber\\
&\P(b \otimes c) = q^{-1/2} c \otimes b,\quad& \P(b \otimes d) = q^{1/2} d \otimes b\nonumber\\
&\P(c \otimes a) = q^{1/2} a \otimes c,\quad& \P(c \otimes b) = q^{-1/2} b \otimes c + q^{-1/2} (q - q^{-1}) a \otimes d\nonumber\\
&\P(c \otimes c) = q^{1/2} c \otimes c,\quad& \P(c \otimes d) = q^{-1/2} d \otimes c + q^{-1/2} (q - q^{-1}) c \otimes d\nonumber\\
&\P(d \otimes a) = q^{-1/2}  a \otimes d,\quad& \P(d \otimes b) = q^{1/2} b \otimes d\nonumber\\
\label{slq2_braiding}
&\P(d \otimes c) = q^{-1/2}  c \otimes d,\quad& \P(d \otimes d) = q^{1/2} d \otimes d
\end{eqnarray}

As a special case of \cite{hayashi}, it is straightforward to show  that:

\begin{prop} For this braiding there are precisely two ribbbon automorphisms $\s_\pm$, defined on generators by $\s_\pm (x) = \pm q^{3/2} x$ for $x = a,b,c,d$, and extended by $\s_\pm (xy) = {\mu} ( \s_\pm \otimes \s_\pm) \P^2 (x \otimes y)$.
\end{prop}

Since in the classical limit $q=1$ we would like $\s = \id$, we will restrict attention to $\s_+$.
The corresponding $\ss \in \A^\circ$ is defined 
 by  $\ss(a) = q^{3/2} = \ss(d)$, $\ss(b) =0 = \ss(c)$ (see \cite{kassel}, p366). 
Then:

\begin{prop} For the braiding (\ref{slq2_braiding}) and $\s = \s_+$, 
$HH^{\P,\s}_0 (\A) = 0$.
\end{prop}
\begin{pf} We calculate $H^{\P,\s}_0 (\A,\A)$ directly from the definition. Since $\A$ is unital this coincides with $HH^{\P,\s}_0 (\A)$. 
 By induction we obtain the formulae
\begin{eqnarray}
&&\P( a^i b^j c^k \otimes  c) = q^{(i-j+k)/2} c \otimes a^i b^j c^k ,\quad
\P( a^i b^j c^k \otimes  a) = q^{(i-j+k)/2} a \otimes a^i b^j c^k,\nonumber\\
&&\P( d^i b^j c^k \otimes  c) = q^{(-i-j+k)/2} c \otimes d^i b^j c^k.\nonumber
\end{eqnarray}
Hence
\begin{eqnarray}
&&b_{\P,\s} ( a^i b^j c^k \otimes a ) = q^{-j-k} ( 1 - q^{(3+i+j+3k)/2} ) a^{i+1} b^j c^k,\nonumber\\
&&b_{\P,\s} ( b^j c^k \otimes c ) = (1 - q^{(3-j+k)/2} ) b^j c^{k+1},\nonumber\\
&&b_{\P,\s} ( d b^j c^k \otimes a) = q^{-j-k-1} [ \; (1 - q^{(j+3k+5)/2} ) + q^{-1} ( 1 - q^{(j+3k+9)/2} ) \; bc \; ] \; b^j c^k\nonumber\\
&&b_{\P,\s} ( b \otimes b) = q^{1/2} b^2, \quad 
b_{\P,\s} [ d \otimes a - q^{-1} ( 1 + q +q^2 ) b \otimes c ] =1.
\nonumber
\end{eqnarray}
Using these in order, first
 $[ a^{i+1} b^j c^k ] = 0$ for all $i,j,k \geq 0$. Second, $[ b^j c^{k+1} ] = 0$ unless $j=k+3$. Next, for all $j,k$, 
 $[ b^j c^k ] = \l_{jk} [ b^{j+1} c^{k+1}]$ for some nonzero $\l_{jk}$, hence  each $[ b^{k+3} c^{k+1}]$ is proportional to $[b^2]$, which is zero.  Finally $[1] =0$, so we have $[ a^i b^j c^k ] =0$ for all $i, j,k \geq 0$. In the same way,
 \begin{eqnarray}
 && b_{\P,\s} ( d^i b^j c^k \otimes c ) = (1 - q^{(3+i-j+k)/2} ) d^i b^j c^{k+1},\nonumber\\
 && b_{\P,\s} ( d^{i+1} b^j c^k \otimes a ) = q^{-j-k} d^{i+1} \; [ \; (1 - q^{(2-i+j+3k)/2} ) + q^{-1} ( 1 - q^{(3i +j + 3k +6)/2} ) \; bc \; ] \; b^j c^k, \nonumber\\
 && b_{\P,\s} ( d^{i+1} b^j \otimes b) = (1 - q^{(3i+j+6)/2} ) d^{i+1} b^{j+1},\quad
 b_{\P,\s} ( d^i \otimes d) = (1 - q^{(i+3)/2} ) d^{i+1}.
 \nonumber
 \end{eqnarray}
 Hence $[d^i b^j c^{k+1} ]=0$ unless $j =i+k+3$. In this case, $[ d^{i+1} b^j c^{k+1} ]$ is proportional to $[ d^{i+1} b^{i+3}]$, which is zero. Finally $[ d^{i+1} b^{j+1} ]=0$ for all $i, j$.
\end{pf}

We now pass to the braided Hopf algebra $\B = \B(\A)$. 
We define new generators
$$u = d,\quad x = q b,\quad y=q c, \quad z =\frac{qa-qd}{q+q^{-1}}$$
Using (\ref{B(A)_product}, \ref{B(A)_antipode},  \ref{B(A)_braiding}) we have  the braided Hopf algebra structure (we drop the ``*" notation for the product)
\begin{eqnarray}
&& ux = q^2 xu, \quad uy = q^{-2} yu, \quad xy = u^2  + (1 + q^{-2} ) uz - 1,\quad zu = uz,\nonumber\\
&& yx = u^2 + (1 + q^2) uz -1, \quad
zx = xz + (1- q^2 ) xu, \quad zy = yz + (1 - q^{-2}) yu\nonumber\\
&& \D(u) = u \otimes u + q^{-2} y \otimes x, \quad
 \D(x) = x \otimes u + u \otimes x + ( 1 + q^{-2}) z \otimes x\nonumber\\
&& \D(y) = y \otimes u + u \otimes y + (1 + q^{-2}) y \otimes z \nonumber\\
&& \D(z) = z \otimes u + u \otimes z + (1 + q^{-2}) z \otimes z + (1 + q^{-2})^{-1} [ x \otimes y -  y \otimes x ]\nonumber\\
&&\underline{S}(u) = u + (1 + q^2 )z,\quad
\underline{S} (x) = - q^2 x, \quad 
\underline{S} (y) = - q^2 y, \quad 
\underline{S} (z) = - q^2 z\nonumber\\
\label{defn_B}
&&\e (u) = 1,\quad \e(x) = \e(y) = \e(z) =0
\end{eqnarray}
Further, we note that $t := u+z = \frac{q a + q^{-1} d}{q+q^{-1}}$ is a central element. 
$\B$ is $ \mathbb{Z} $-graded with 
$x$, $y$ having degree 1, -1, and $u,z$ having degree zero. 
Using this and the commutation relations gives that 
$$ e_{ijk} := \left\{\begin{array}{cc}
				x^i u^j z^k & : i \geq 0 \cr
				y^i u^j z^k & : i \leq 0 \end{array}\right. \quad j, k \in \bN  $$
is a vector space basis of $\B$. 
The braiding (\ref{B(A)_braiding}), which for convenience we denote by $\P$ rather than $\P_B$, is given by: 
\begin{eqnarray}
&& \P ( x \otimes x )= q^2 x \otimes x, \quad
 \P (x \otimes y )= q^{-2} y \otimes x, \quad
 \P (x \otimes z )= z \otimes x, \quad
 \P (x \otimes u )= u \otimes x\nonumber\\
&& \P (y \otimes x )= q^{-2} x \otimes y + (1- q^{-2}) f(2) \;  y \otimes x  - (1 + q^{-2}) f(2) \; z \otimes z\nonumber\\
&& \P (y \otimes y )= q^2 y \otimes y, \quad
 \P (y \otimes z )= z \otimes y + f(2) y \otimes z, \quad
 \P (y \otimes u )= u \otimes y - f(2) y \otimes z\nonumber\\
&& \P (z \otimes x )= x \otimes z + f(2) z \otimes x,\quad
\P (z \otimes y )= y \otimes z \nonumber\\
&& \P (z \otimes z )= z \otimes z + (q^{-2}-1) y \otimes x, \quad  \P (z \otimes u )= u \otimes z + (1 - q^{-2}) y \otimes x\nonumber\\
&& \P (u \otimes x )= x \otimes u - f(2) z \otimes x , \quad
 \P (u \otimes y )= y \otimes u, \quad
\P (u \otimes z )= z \otimes u + (1- q^{-2}) y \otimes x, \nonumber\\
\label{braided_slq2_braiding}
&&\P (u \otimes u )= u \otimes u + ( q^{-2} -1 ) y \otimes x
\end{eqnarray}
where $f(n) = q^n - q^{-n}$.
It follows that $(\id \otimes \underline{S}) \P = \P (\underline{S}  \otimes \id)$,  $(\underline{S}  \otimes \id) \P = \P (\id \otimes \underline{S})$.

\begin{lemma} There  are precisely two  ribbon automorphisms $\s_\pm$ of $\B$, given by
\begin{equation}\label{sigma_pm}
 \s_\pm (u) = \pm [ u + (1 - q^4) z ], \quad \s_\pm (x) = \pm q^4 x, \quad \s_\pm (y) = \pm q^4 y, \quad \s_\pm (z) = \pm q^4 z
\end{equation} 
(this implies $\s_\pm (t) = \pm t$)
equivalently by
$$\s_\pm (a) = \pm [ q^2 a + (1- q^2 ) d], \quad \s_\pm (b) = \pm q^4 b, \quad \s_\pm (c) = \pm q^4 c, \quad \s_\pm (d) = \pm [ (q^2 - q^4) a + {\frac{(q^6 +1)}{q^2 +1}} d] $$
\end{lemma}
\begin{pf} Demanding compatibility of $\s$ with the defining relations, for example 
$$0 = \s(xy - yx + (q^2 - q^{-2}) uz ) = {\mu}  (\s \otimes \s) [ x \otimes y - y \otimes x + (q^2 - q^{-2}) z \otimes u ]$$
gives $\s(x) = \l_1 x, \quad \s(y) = \l_2 y,\quad  \s(z) = \e q^4 z, \quad \s(u) = \e[ u + (1 - q^4 ) z]$, where $\e = \pm 1$ and $\l_1 \l_2 = q^8$.
 Compatibility with the braiding forces $\l_1 = \l_2 = \e q^4$, hence the result. 
\end{pf}

It is natural to require that the ribbon automorphism becomes the identity in the classical limit $q=1$, which imposes $\e =1$. Hence we work with  the ribbon automorphism $\s = \s_+$.
 The ground field $\bC$ becomes a left $\B$-module through the
character $ \varepsilon $. It is easy to check that:

\begin{lemma} The following is a resolution of $\bC$ by free left $\B$-modules:
\begin{equation}\label{free_B_res} 
     	\begin{CD}
	{0} @ >>> {B} @ >{\varphi_3} >>{B^3} @ >{\varphi_2} >>{B^3} @ >{\varphi_1} >> {B} @ >{\varphi_0} >> {\field} @ >>> {0} @. .\\
    	\end{CD}
\end{equation}
where $\varphi_0 (a) = \e(a)$, $\varphi_1 (a,b,c) = ax + by + c(u-1)$, $\varphi_3(a) = a( y, -q^2 x, u-1)$ and 
$$\varphi_2 (a,b,c) = (a,b,c) 
\left( \begin{array}{ccc}
q^{-2} u & 0 & -x\cr
0 & q^2 u -1 & -y\cr
-y & q^2 x & (1- q^2) (u+1)\cr
\end{array}\right)$$
\end{lemma}

Given the resolution, we can compute braided Hochschild homology, giving a ``no dimension drop" result along the lines of \cite{bz,tompodles,hk,hk2}:

\begin{thm}
\label{HH3B_SLq2}
 For the braiding  (\ref{braided_slq2_braiding}) and $\s = \s_+$ (\ref{sigma_pm}), 
$H_3^{\P,\s} (\B,\B) \cong \bC$.
\end{thm}
\begin{pf}
 Tensoring (\ref{free_B_res}) on the left by $R(B) \otimes_B -$ gives 
 $$H_3^{\P,\s} (\B,\B) = \ker \{ \id \otimes \varphi_3 :  R(B) \otimes_B B \to R(B) \otimes_B B^3 \} = \ker \{ {\tilde{\varphi}} : B \to  B^3 \}$$
 where $\tilde{\varphi} (a) = ( a \btl y, - q^2 a \btl x, a \btl (u-1))$.
  To compute the right action $\btl$  (\ref{black_right_action_A_on_M}) of $x$, $y$, $u$ on PBW monomials $a = x^i u^j z^k$, $y^i u^j z^k$, we need to compute the braidings $\P ( a \otimes t)$ for $t = x,y,u$. 
 Lengthy but straightforward calculations give  the formulae:
 \begin{eqnarray}
 &\P(x^i u^j z^k \otimes y) =& q^{-2i} y \otimes x^i u^j z^k, \quad \P(y^i u^j z^k 
 \otimes  y) = q^{2i} y \otimes y^i u^j z^k,\nonumber\\
 &\P( x^i u^j z^k \otimes u) =& u \otimes x^i u^j z^k  + q^{-2i} y \otimes x^i [ q^{-2} x u^{j-1} z^k - u^{j-1} z^k x]\nonumber\\
 &\P(y^i u^j z^k \otimes u) =& u \otimes y^i u^j z^k  + (1 - q^{2i}) y \otimes y^{i-1} u^j z^k [ (1 + q^{-2}) z + ( q^{-2} - q^{-2i}) u ] \nonumber\\
 &&  \quad + q^{2i}  y \otimes y^i u^{-1} [ x u^j z^k - u^j z^k x ] \nonumber
 \end{eqnarray}
 where $u^{-1} [ x u^j z^k - u^j z^k x ]$ is notational shorthand for 
 $$q^{-2} x u^{j-1} z^k - u^{j-1} z^k x = q^{-2} x u^{j-1} [ z^k - q^{2j} ( z + (1- q^2 ) u )^k ]$$
 which is well-defined even for $j=0$, being in this case equal to 
 $- q^{-2} \sum_{l=1}^k {\tiny\kchoosel} (1-q^2 )^{k-l} x z^l u^{k-l-1}$, with the empty sum ($k=0$) being taken to be zero. 
 Furthermore,
$$\s \underline{S}^{-1} (u) =  u + (1+ q^2) z,\quad  
\s \underline{S}^{-1} (x) =  - q^2 x,\quad  
\s \underline{S}^{-1} (y) =  - q^2 y,\quad
\s \underline{S}^{-1} (z) =  - q^2 z$$
 It follows that the actions of $u$ and $y$ on PBW monomials are:
\begin{eqnarray}
&& x^i u^j z^k \btl u = q^{-2i} x^i u^j z^k, \quad y^i u^j z^k \btl u = q^{2i} y^i u^j z^k\nonumber\\
&& y^i u^j z^k \btl y = (q^2 - q^{2i}) (1 - q^{-4i}) y^i u^{j+1} z^k y + q^{2i} y^i u^{-1} [ u^j z^k y - y u^j z^k ]\nonumber
\end{eqnarray}
Therefore  $a \btl (u-1) = 0$ if and only if $a = \sum \a_{j,k} u^j z^k$, for some $\a_{j,k} \in \bC$. Now, 
$$u^j z^k \btl y = u^{j-1} z^k y - q^2 y u^{j-1} z^k = q^2 y u^{j-1} [ q^{-2j} ( z + (1- q^{-2} ) u )^k - z^k ]$$
Hence $(  \sum \a_{j,k} u^j z^k ) \btl y =0$ if and only if $\a_{j,k} =0$ for $(j,k) \neq (0,0)$. 
So $a \btl (u-1) = 0 = a \btl y$ if and only if $a = \l 1$.
Finally it is easy to check that $1 \btl x =0$. 
Hence $\ker ( \tilde{\varphi} ) = \bC [1]$. 
\end{pf}

\section{Acknowledgements}

We both thank Shahn Majid for very useful discussions,
and the Isaac Newton Institute, Cambridge for hosting
us during the time this work was completed. We also
thank the referee for their careful reading of the paper
and many useful comments and suggestions.
T.H.~thanks the EPSRC for their support via a Postdoctoral Fellowship,  Katedra Metod Matematycznych Fizyki, Uniwersytet Warszawski for its support via EU Transfer of Knowledge contract MKTD-CT-2004-509794, and Instytut Matematyczny, Polska Akademia Nauk   for their hospitality.
He also thanks the School of Mathematical Sciences, Queen Mary, University of London for their support. 
 U.K.~thanks the EU for support via Marie Curie EIF
 515-144 and the EPSRC for support via EP/E/043267/1.

\end{document}